\documentclass[leqno,11pt,a4paper]{amsart}
\usepackage{tgbonum}
\usepackage{amsmath}
\usepackage{amsthm}
\usepackage{amsfonts,amssymb}
\usepackage{bbm}
\usepackage{multirow}
\usepackage{latexsym}
\usepackage{mathrsfs}
\usepackage{stackengine}
\usepackage[all]{xy}
\usepackage{url}
\usepackage{xcolor}
\usepackage{eqparbox}
\usepackage{tikz}
\usepackage{float}
\usepackage{hyperref}
\usepackage{amssymb,longtable}
\usepackage{fontenc}
\usepackage{array,etoolbox}
\preto\tabular{\setcounter{magicrownumbers}{0}}
\newcounter{magicrownumbers}

\usepackage{pdfsync}

\usepackage{amsthm}

\usepackage{hyperref}
\usepackage{enumerate}
\usepackage{url}

\usepackage{longtable}
\usepackage{booktabs}
\usepackage{colortbl}
\newcommand{\evnrow}{\rowcolor[gray]{0.95}}
\newcommand{\oddrow}{}

\usepackage{geometry}

\usepackage{pdflscape}
\usepackage{longtable}
\usepackage{booktabs}
\usepackage{color}
\usepackage{arydshln}
\usepackage{verbatim}
\usepackage{enumitem}\usepackage{mathtools}

\usepackage{mathdots}
\usepackage[english]{babel}
\usepackage{lipsum}
\usepackage{amsthm}
\usepackage{todonotes}
\usepackage[vcentermath]{youngtab}

\newcommand\stackrqarrow[2]{%
    \mathrel{\stackunder[2pt]{\stackon[4pt]{$\rightsquigarrow$}{$\scriptscriptstyle#1$}}{%
            $\scriptscriptstyle#2$}}}

\newcommand{\eqmathbox}[2][M]{\eqmakebox[#1]{$\displaystyle#2$}}

\theoremstyle{plain}
\newtheorem{lema}{Lemma}[section]

\newtheorem{thrm}[lema]{Theorem}

\theoremstyle{remark}

\newtheorem{rmk}[lema]{Remark}

\theoremstyle{definition}
\newtheorem{dfn}[lema]{Definition}

\def\bQ{{\mathbb Q}}
\def\bC{{\mathbb C}}
\def\bN{{\mathbb N}}

\def\bP{{\mathbb P}}
\def\bV{{\mathbb V}}
\def\bI{{\mathbb I}}

\def\QQ{{\mathbb Q}}
\def\CC{{\mathbb C}}

\def\AA{{\mathbb A}}

\def\PP{{\mathbb P}}
\def\VV{{\mathbb V}}
\def\II{{\mathbb I}}

\def\p2xp2{\mathbb{P}^2\times\mathbb{P}^2}
\def\pxp{\mathbb{P}^2\times\mathbb{P}^2}

\def\wp2{\mathrm{w}($\bP^2\times \bP^2$)}

\def\Oh{\mathcal{O}}

\def\gen{{\mathrm{gen}}}
\def\wt{{\mathrm{wt}}}

\def\w{{\mathrm{w}}}

\def\gen{{\mathrm{gen}}}

\def\Y{{\mathbf Y}}

\def\Tom{\mathrm{Tom}}
\def\Jer{\mathrm{Jer}}

\def\Pf{{\mathrm{Pf}}}

\makeatletter
\newcommand*{\dashdownarrow}{%
  \mathrel{%
    \mathpalette\dasharrow@vert{-90}%
  }%
}
\newcommand*{\dashuparrow}{%
  \mathrel{%
    \mathpalette\dasharrow@vert{90}%
  }%
}
\newcommand*{\dasharrow@vert}[2]{%
  \sbox0{$#1\vcenter{}$}%
  \sbox2{$#1\dashrightarrow\m@th$}%
  \dimen@=1.2\dimexpr\ht2-\ht0\relax
  \sbox2{\raisebox{-\ht0}{\unhcopy2}}%
  \ht2=\z@
  \dp2=\z@
  \vcenter{\hbox to 2\dimen@{\hfill\rotatebox{#2}{\box2}\hfill}}%
}
\makeatother

 \newcommand{\into}{\hookrightarrow}

 \newcommand{\cO}{\mathcal O}
 \newcommand{\cC}{\mathcal C}
 
 \newcommand{\cF}{\mathcal F}

 \newcommand{\Ybar}{\overline Y}
 
 \newcommand{\Ygen}{ Y_\gen}

 \newcommand{\si}{\sigma}

 \newcommand{\al}{\alpha}

\newcommand{\PxP}{\bP^2 \times \bP^2}

\title{Construction and deformations of Calabi--Yau 3-folds in  codimension 4 }
\author[S.~Mohsin, S.~Nazir, M.I.~Qureshi]{Sumayya Mohsin, Shaheen Nazir, and Muhammad Imran Qureshi}
\address{Sumayya Mohsin\newline
Department of Mathematics, 
Lahore University of Management Sciences (LUMS), Lahore, Pakistan}
\email{17070002@lums.edu.pk}
\address{Shaheen Nazir\newline
Department of Mathematics, 
Lahore University of Management Sciences (LUMS), Lahore, Pakistan}
\email{shaheen.nazir@lums.edu.pk}
\address{Muhammad Imran Qureshi\newline
Department of Mathematics, King Fahd University of Petroleum \& Minerals (KFUPM), Dhahran
31261, Saudi Arabia \&\newline
Interdisciplinary Research Center for Intelligent and Secure Systems, King Fahd University of Petroleum \& Minerals (KFUPM), Dhahran
31261, Saudi Arabia}
\email{imran.qureshi@kfupm.edu.sa}
\begin{document}

\maketitle  
\begin{abstract} We construct polarized Calabi--Yau 3-folds  with at worst isolated canonical orbifold points in codimension 4 that can be described in terms of the equations of the Segre embedding of $\pxp$ in $\PP^8$.  
We investigate the existence of other deformation families  in their Hilbert scheme   by either studying Tom and Jerry degenerations or by comparing their Hilbert series with  those of existing low codimension Calabi--Yau 3-folds.  Among other interesting results, we find a  family  of Calabi--Yau 3-fold with five distinct Tom and Jerry deformation families, a phenomenon not seen for $\QQ$-Fano 3-folds. We compute the Hodge numbers of     \(\p2xp2 \) Calabi--Yau 3-folds and corresponding manifolds obtained by performing  crepant resolutions. We obtain a manifold with   a   pair of Hodge numbers that does  not
 appear in the  famously known list of 
30108 distinct Hodge pairs  of Kruzer--Skarke, in  the list of 7890 distinct Hodge pairs
corresponding to    complete intersections in the product of projective spaces
and in Hodge paris obtained from Calabi--Yau 3-folds having  low codimension embeddings in  weighted projective spaces.           \end{abstract}


\section{Introduction}

\subsection{Calabi--Yau 3-folds and Gorenstein rings}

 In this article, a Calabi--Yau (CY) 3-fold  is a three dimensional normal projective algebraic variety $X$  such that \(H^1(X,\cO_X)=0\) and the  canonical divisor class \(K_X\cong \mathcal{O}_X \). We allow $X$ to contain at worst canonical isolated orbifold points and so it admits a crepant resolution to a smooth CY 3-fold. We use $X$ for a family of Calabi--Yau 3-folds and one specific representative of the family interchangeably, if no confusion can arise.

Over the last 4 decades, the study of CY 3-folds has been an important area of research in algebraic geometry and string theory, paving the way for  connections between the two areas. A key role of CY 3-folds in string theory was first highlighted by Candelas-Horowitz-Strominger-Witten in \cite{Candelas-1}. Since then the construction  and  classification of CY 3-folds have been a subject of unifying interest in both areas of research. Some of the initial lists  of        CY 3-folds were constructed  as complete intersections in the product of projective spaces by     Candelas--Dale--L\"utken--Schimmrigk \cite{Candelas-CICY} and by  Candelas-Lynker-Schimmrigk  \cite{Candelas-wp4} as hypersurfaces in the weighted projective space \(\w\PP^4\).

We construct Calabi--Yau 3-folds as polarized orbifolds \((X, D)\): \(X\) is a normal projective algebraic variety with a choice of an ample \(\bQ\)-Cartier divisor \(D\). A choice of \(D\) gives rise to a graded ring \[R(X,D)=\bigoplus_{m\in \bN} H^0\left(X, \Oh_X(mD)\right).\]
A surjective morphism \[k[x_1,\ldots,x_n]\twoheadrightarrow R(X,D)\]
from a free graded polynomial ring with degree\((x_i)=a_i\), gives the embedding of \(X\) in the  weighted projective space \(\bP(a_1,\ldots,a_n). \)  If  the image of \(X\) is described by a minimal generating set  \(f_1,\ldots,f_k\) and \(\II(X)=(f_1,\ldots,f_k)\)  then  \[R(X,D)\cong \dfrac{k[x_1,\ldots,x_n]}{\II(X)}. \]
We are interested in the case when \(X\) is  a wellformed,  quasismooth,
and  projectively Gorenstein   variety. Recall that a weighted projective variety $X\subset \w\bP$ of codimension $c$ is  wellformed if it
does not contain a singular locus of dimension greater than or equal to $c - 1$, 
it is  quasismooth if the punctured affine cone \(\widetilde X\backslash\{0\}\)
contains no singularities, and it is  projectively Gorenstein if   \(R(X, D)\)  is Cohen--Macaulay and the canonical divisor  \(K_X\) is subcanonical.

 The case of CY 3-fold  hypersurfaces ($\II(X)=(f))$ in \(\w\bP^4\) was initially studied by Candelas et al in \cite{Candelas-wp4} and completely classified by Kruezer--Skarke in \cite{KS}, resulting in  7555 orbifolds \cite{GRDB}.   In codimension 2, Gorenstein varieties  are complete intersections,  i.e. defined by two equations; and in codimension 3 they are described by the maximal Pfaffians of a \((2n+1)\times (2n+1) \) skew-symmetric matrix by  Buchsbaum--Eisenbud  \cite{BE}. These algebraic structures  have been utilized to construct many models of  CY 3-folds $\w\PP^6$ in \cite{tonoli,kapustka,kapustka2,BKZ} and elsewhere. 

There is currently no ``user-friendly'' structure theorem in codimension greater than or equal to 4, though progress on the  codimension 4 case has been made by  Reid \cite{codim4}. Therefore a variety of ideas have been employed to  construct and classify smooth codimension 4  Gorenstein  CY 3-folds  in  \(\PP^7\) by Bertin  \cite{Bertin}, Caughlan--Golebiowski--Kapustka--Kapustka \cite{CGKK} and    Schenck--Stillman--Yuan \cite{SSY}; followed by subsequent codimension 4 Gorenstein constructions in \cite{kapustka2021quaternary, ablett2021halfcanonical, ablett2024deformations}.  In the weighted $\w\PP^7$  case, some  lists of Gorenstein   polarized CY 3-folds in codimension greater than or equal to 4 have been constructed by using techniques of  Gorenstein formats \cite{QS,QS-AHEP,QS2,BKZ}   and parallel type-I unprojection  \cite{BG}.                                

  We focus on the construction  of  deformation families of Gorenstein canonical CY 3-folds  in codimension 4 
 such that one of the  deformation families  \(X\)   is described by  using equations of the Segre embedding of \(\pxp\), i.e. as \(2\times 2\) minors of the size 3 generic square matrix. Equivalently, we can consider \(X\) as a (weighted) complete intersection of some projective cone(s) over a weighted \(\pxp\) variety.      We say that \(X\) is a \(\pxp\) Calabi--Yau 3-fold or \(X\) lives in \(\pxp\) Gorenstein format~\cite{BKQ}.   
\subsection{\(\pxp\) Calabi--Yau 3-folds}
We first recall the definition of weighted  \(\pxp\) variety. 
\begin{dfn}
Let  $u=(u_1,u_2,u_3)$ and $v=(v_1,v_2,v_3)$ be two vectors of integers or half integers satisfying 
$$u_1+v_1\geq 1,\ \ \  u_1\leq u_2\leq u_3,\ \ \ \ v_1\leq v_2\leq v_3. $$
 Consider the Segre embedding $$j:\p2xp2:=W \into\PP^8(x_{ij}:1\le i,j\le 3)$$   then the image of the 4-fold
 embedding  can be  described by  
$$
\Bigg(\bigwedge^2\begin{pmatrix}
x_{11} &x_{12} &x_{13}\\
x_{21} &x_{22} &x_{23}\\
x_{31} &x_{32} &x_{33}\\
\end{pmatrix}=0\Bigg)\subset\bP^8(x_{ij}).
    $$
The \textit{weighted \(\pxp\) variety} \(\cF\) is the  quotient of the punctured  affine cone $\widetilde{W}\setminus \{0\}  $ by the $\bC^{\times}$ action  : $$\epsilon:x_{ij}\rightarrow \epsilon^{a_{ij}}x_{ij}, $$ where \(a_{ij}=u_i+v_j\) for all $1\leq i,j\leq 3 $, giving the embedding \(\cF\into\PP^8(a_{ij})\).
If \(\cF\) is wellformed then the canonical divisor class  \(K_X=\Oh\left(-\displaystyle\sum_{j=1}^3
a_{jj}\right)\).
\end{dfn}

  We start by  producing a list of  candidate families of isolated canonical CY 3-folds in \(\pxp\) format by using an algorithm developed in \cite{QJSC, BKZ}. In practice, the algorithm searches for the Hilbert  series corresponding  to some equations format (like minors of a matrix, Pfaffians of a matrix, etc.) that may be the Hilbert series of a polarized orbifold with isolated orbifold points. We apply it to search for those Hilbert series in codimension four  \(\pxp\) format, that may be realized by some isolated canonical  CY 3-fold.

  We focus on a detailed study of CY 3-folds with relatively small sets of weights of \(\bP(a_1,\ldots,a_8)\). In particular, we restrict our attention to  CY 3-folds \(X\) in \( \bP(a_1,\ldots,a_8)\) with \(\sum a_i\le 20\) and  we get 24 Hilbert series that satisfy the above condition. We   carry out a complete singularity  and quasismoothness analysis of the corresponding equations obtained from \(\pxp\) format.  Among these, 23 of these give rise  to    wellformed and quasismooth \(\pxp\) CY 3-folds with at worst isolated canonical orbifold points, and one of the candidates contains worse singularities.
\begin{thrm}
     \label{Th!pxp}Let \(X\) be a family of polarized Calabi--Yau 3-folds with at worst isolated canonical orbifold points embedded as codimension 4 variety in \(\bP(a_1,\ldots,a_8)\) in \(\pxp\) Gorenstein format. Then there exist  at least 23 such families of  wellformed and quasismooth Calabi--Yau 3-folds in \(\bP(a_1,\ldots,a_8)\).  
 
\end{thrm}
Out of 23 working cases, 4 of them are smooth CY\ 3-folds, i.e. three dimensional Calabi--Yau  manifolds. The rest of the 19 families are orbifolds and they admit a crepant resolution to a smooth CY manifold. 
 We then compute the Hodge numbers of each $\pxp$  CY 3-fold \(X\) and its crepant resolution \(\hat X\) by using computer algebra, as explained in Section \ref{S!Hodge}. The  resolution \(\hat X\) of the family no. 7 (Table \ref{T!Summary}) is a  smooth CY 3-fold with Hodge numbers \(h^{1,1}(\hat X)=3\) and \(h^{2,1}(\hat X)=62\). The Hodge pair \((3,62)\)  appears neither in the Kruzer--Skarke \cite{KS} list of distinct  30108 Hodge pairs nor in the list of complete intersections in the product of projective spaces  list of 7890 distinct Hodge pairs in \cite{CICY-List}, giving it a distinctive position in   the Calabi--Yau landscape~\cite{He-Book}.
 
\subsection{Deformation families   of \(\pxp\) CY 3-folds} 
 We  investigate the existence of other deformation families of every \(\pxp\) CY 3-fold \(X\) of Theorem \ref{Th!pxp} by using the following two ideas. 
 \begin{enumerate}
\item According to \cite{BKR}, a  typical codimension four  \(\pxp\) variety with   a numerical type-I center, comes from a special Tom unprojection of a codimension 3 Pfaffian variety. 
So our first  approach is to construct other  possible Tom and Jerry families of \(X\);  requiring analysis of the corresponding  Pfaffian  CY 3-fold in codimension 3 that is obtained from the projection of type-I center of $\pxp$ CY 3-fold.
  \item  The other possibility is to  search for low codimension models by matching the Hilbert series of \(\pxp\) models with  existing classes of low codimension CY 3-folds and show that  they lie in different deformation families. \end{enumerate} 


We start our deformation analysis with  a quasismooth and wellformed \(\p2xp2\) CY 3-fold \(X\) with isolated canonical orbifold points embedded in the weighted projective space $\bP(a_1,\ldots,a_7,r )$. Recall that  a point \(p:=\frac 1r(a,b,c)\) is numerical type-I center if the  local variables \(x_1,x_2\) and \(x_3\)  of weight $a,b $ and $c$ respectively   are part of basis of vector spaces \(H^0(aD), H^0(bD)\) and \(H^0(cD)\),
 for some $\{a,b,c\}\subset \{a_1,\ldots,a_7\}$. If \(X\) contains a numerical type-I center $p$,  then we perform a Gorenstein projection  to a   codimension 3  CY 3-fold  \[(\overline{Y} \supset D)\into \bP(a_1,\ldots,a_7)\] which maps \(p\) to the complete intersection plane divisor \(D:=\bP(a,b,c)\) that is projectively normal in \(\Ybar\).  The equations of \(\Ybar\) can be written as the maximal Pfaffians of a \(5\times 5\)  skew-symmetric matrix \(M\),   in   \(\Tom_i\) format for some \(i \in \{1,\ldots,5\}\). If \(\Ybar\) contains at worst nodes as its singularities then we say that \(X\) is a CY 3-fold of \(\Tom_i\) type.    This is followed by a  deformation \(\Ybar \to Y\) of \(\Tom_i\) to a Pfaffian CY 3-fold   \( Y\) by deforming the    entries of \(M\) to generic forms in the graded polynomial ring corresponding to \(\bP(a_1,\ldots,a_7)\). Then we  perform all possible Tom and Jerry   degenerations of \(Y\) and all distinct nodal CY 3-folds  give an unprojection to a polarized CY 3-fold in codimension 4 having the same  numerical data (Hilbert series, Hilbert polynomial, orbifold points) but different Euler characteristics, i.e. they lie in  different components of the Hilbert scheme of \(\p2xp2\) CY 3-folds. 

Our analysis partly relies on the classical method of Hirzerbruch~\cite{Hirz}:  constructing nodal CY 3-folds and then performing resolution of these singularities to get new  CY 3-folds,   used  in various contexts,
see
  \cite{BG, CY3-nodes} for example. 

The following table summarizes the results of our analysis, including the number of deformation families.

\begin{longtable}{>{\hspace{0.5em}}llllccc<{\hspace{0.5em}}}
\caption{Summary of  results} \label{T!Summary}\\
\hline & \\[-1.5ex]
HS & Embedding & Orbifold points & $D^3$ & $h^{2,1}(X)$ & $h^{1,1}(\hat X)$  &$\#$Fam.\\
\cmidrule(lr){1-1}\cmidrule(lr){2-4}\cmidrule(lr){5-6}\cmidrule(lr){7-7} 
\endfirsthead
\hline & \\[-1.5ex]
HS & Embedding & Orbifold points & $D^3$ & $h^{2,1}(X)$ & $h^{1,1}(\hat X)$  &$\#$Fam.\\
\cmidrule(lr){1-1}\cmidrule(lr){2-4}\cmidrule(lr){5-6}\cmidrule(lr){7-7} \endhead
\multicolumn{7}{r}{\scriptsize\emph{Continued on next page}}\\
\endfoot
\bottomrule
\endlastfoot
\evnrow 1 & 
   $ \begin{array}{@{}l@{}} X_{2^3,3^6}\\\quad\subset\bP^7\end{array}$ & & $17$ &  $58$ &\(2\)& $1$\\ 

 \oddrow  2 &  $ \begin{array}{@{}l@{}} X_{2^3,4^6}\\\quad\subset\bP(1^7,3)\end{array}$ & $\frac{1}{3}(1,1,1)$ & $\frac{34}{3}$& $65$ &$3$& $2$\\

\evnrow 3 &  $ \begin{array}{@{}l@{}} X_{3^6,4^3}\\\quad\subset\bP(1^6,2^2)\end{array}$ & & $11$&$50$ &$2$& $1$\\

\oddrow 4 &  $ \begin{array}{@{}l@{}} X_{2,3^4,4^4}\\\quad\subset\bP(1^7,3)\end{array}$ & $\frac{1}{3}(1,1,1)$ & $\frac{40}{3}$ & $78$ &$3$& $4$\\

\evnrow 5 &$ \begin{array}{@{}l@{}} X_{2,3^4,4^4}\\\quad\subset\bP(1^6,2^2)\end{array}$ &&  $10$&$50$ &$2$& $2$\\

\oddrow 6 &$ \begin{array}{@{}l@{}} X_{4^9}\\\quad\subset\bP(1^4,2^4)\end{array}$ & & $6$  &$38$&$2$& $1$\\

\evnrow 7 & $\begin{array}{@{}l@{}} X_{3^2,4^5,5^2}\\\quad\subset\bP(1^5,2^2,3)\end{array}$ & $\frac{1}{3}(1,1,1)$ & $\frac{22}{3}$& $62$ &$3$& $5$\\

\oddrow 8 & $ \begin{array}{@{}l@{}} X_{4^6,6^3}\\\quad\subset\bP(1^5,3^3)\end{array}$ & $3\times\frac{1}{3}(1,1,1)$ & $6$ & $62$ &$5$& $3$\\

\evnrow 9 & $ \begin{array}{@{}l@{}} X_{2,3^2,5^2,6^4}\\\quad\subset\bP(1^5,2^2,5)\end{array}$ & $\frac{1}{5}(1,2,2)$ & $\frac{29}{5}$& $84$ &$4$& $1$\\

\oddrow 10& $ \begin{array}{@{}l@{}} X_{3,4^3,5^3,6^2}\\\quad\subset\bP(1^4,2^2,3^2)\end{array}$ & $\frac{1}{3}(1,1,1)$  & $\frac{13}{3}$& $61$ &$3$& $3$\\

\evnrow 11& $ \begin{array}{@{}l@{}} X_{2,4^4,6^4}\\\quad\subset\bP(1^5,3^3)\end{array}$ & $2\times\frac{1}{3}(1,1,1)$ & $\frac{14}{3}$& $97$ &$4$& $2$\\

\oddrow 12 &  $ \begin{array}{@{}l@{}} X_{2,4^4,6^4}\\\quad\subset\bP(1^6,3,5)\end{array}$ & $\frac{1}{5}(1,1,3)$ & $\frac{42}{5}$& $73$ &$4$& $3$\\

\evnrow 13 &  $ \begin{array}{@{}l@{}} X_{4^4,5^4,6}\\\quad\subset\bP(1^4,2^2,3^2)\end{array}$ & $2\times\frac{1}{3}(1,1,1)$ & $\frac{14}{3}$& $48$ &$4$& $4$\\

\oddrow 14 & $ \begin{array}{@{}l@{}} X_{3,4^2,5,6^3,7^2}\\\quad\subset\bP(1^4,2^2,3,5)\end{array}$ & $\frac{1}{3}(2,2,2)$, $\frac{1}{5}(1,1,3)$ & $\frac{61}{15}$& $67$ &$5$& $1$\\

\evnrow 15 & $ \begin{array}{@{}l@{}} X_{4,5^4,6^4}\\\quad\subset\bP(1^3,2^2,3^3)\end{array}$ & $2\times\frac{1}{3}(1,1,1)$ & $\frac{8}{3}$& $44$ &$4$& $4$\\

\oddrow 16 & $ \begin{array}{@{}l@{}} X_{6^9}\\\quad\subset\bP(1^3,3^5)\end{array}$ & $6\times\frac{1}{3}(1,1,1)$ & $2$& $35$ &$8$& $2$\\

\evnrow 17 & $ \begin{array}{@{}l@{}} X_{2,3^2,7^2,8^4}\\\quad\subset\bP(1^4,2^2,3,7)\end{array}$ & $\frac{1}{3}(1,1,1)$, $\frac{1}{7}(2,2,3)$ & $\frac{58}{21}$& $114$ &$6$& $1$\\

\oddrow 18 & $ \begin{array}{@{}l@{}} X_{4^2,6^5,8^2}\\\quad\subset\bP(1^4,3^3,5)\end{array}$ & $2\times\frac{1}{3}(1,1,1)$, $\frac{1}{5}(1,1,3)$& $\frac{46}{15}$& $73$ &$6$& $4$\\

\evnrow 20 & $ \begin{array}{@{}l@{}} X_{2,5^4,8^4}\\\quad\subset\bP(1^5,2,4,7)\end{array}$ & $\frac{1}{7}(1,2,4)$ & $\frac{36}{7}$& $89$ &$5$& $2$\\

\oddrow 21 & $ \begin{array}{@{}l@{}} X_{4,5^2,6^3,7^2,8}\\\quad\subset\bP(1^3,2^2,3^2,5)\end{array}$ & $\frac{1}{5}(1,1,3)$ & $\frac{12}{5}$& $52$ &$4$& $4$\\

\evnrow 22 & $ \begin{array}{@{}l@{}} X_{4,5^2,6^3,7^2,8}\\\quad\subset\bP(1^3,2,3^3,4)\end{array}$ & $3\times\frac{1}{3}(1,1,1)$ & $2$& $57$ &$5$& $4$\\

 23 & $ \begin{array}{@{}l@{}} X_{3,4^2,7,8^3,9^2}\\\quad\subset\bP(1^3,2^2,3^2,7)\end{array}$ & $\frac{1}{3}(1,1,1)$, $\frac{1}{7}(1,3,3)$ & $\frac{40}{21}$& $86$ &$6$& $1$\\

\evnrow 24 & $ \begin{array}{@{}l@{}} X_{4,5,6^2,7^2,8^2,9}\\\quad\subset\bP(1^3,2,3^2,4,5)\end{array}$ & $\frac{1}{3}(1,1,1)$, $\frac{1}{5}(1,1,3)$& $\frac{26}{15}$& $69$ &$5$& $3$\\
\end{longtable}
  The first column lists the number of the Hilbert series from Table \ref{main
  table} and it does not list 19 since that contains worse singularities. The  data in columns \(2-4\) provides the equations degrees, weights of the ambient weighted projective space, canonical orbifold points, and the self intersection number \(D^{3}\)   for all  codimension 4 deformation families of \(X\). The column \(5-6\) represents the data associated to the \(\pxp\) family \(X\): the Hodge number \(h^{2,1}\) of \(X\),  and Hodge number \(h^{1,1}\) of its crepant resolution \(\hat X\). The Hodge number \(h^{1,1}\) of \(X\) is equal to  2 for all \(\pxp\) families so we do not list it for each case. The last column lists the number of deformation families found for each Hilbert series given in Table \ref{main table}.     

\subsection{Nature of   CY 3-folds}
In  our analysis of studying  deformations  of these  23 families of CY 3-folds, we come across various  interesting phenomena.  Expectedly some of these properties resemble those of \(\bQ\)-Fano 3-folds \cite{BKQ}, but we also
encounter some new phenomena that were not seen in the Fano 3-fold case. We provide a summary of  various novelties,    and results  obtained in our  analysis of these CY 3-folds.
\subsubsection{\textbf{CY 3-folds with five distinct Tom and Jerry families}}
In the study of \(\bQ\)-Fano 3-folds   at most four  distinct Tom and Jerry families were found \cite{BKRbigtable} for each Hilbert series. We prove the existence of a CY 3-fold family with 5 distinct Tom and Jerry deformation families in  Section \ref{S!5TJ}, with the following details.  
\begin{thrm} \label{Th!TJ5}Let \(X\) be a Calabi--Yau 3-fold given with  the Hilbert series \[P_{7}(t)=\frac{1- 2t^3- 5t^4+2t^5+8t^6+2t^7-5t^9-2t^9+t^{12}}{(1-t)^5(1-t^2)^2(1-t^3)}\]  with a single orbifold point \(\frac13(1,1,1)\) and \(D^3=\frac{22}{3}\). Then there exist two  Tom and three Jerry deformation families: realizing  \(P_7(t)\) by five  deformation families, distinguished by their Euler characteristics. 
 $$\begin{tabular}{|c|ccccc|}\hline
\textrm{Family} & $\Tom$ & $\Tom$ & $\Jer$& $\Jer$&$\Jer$\\\hline
\evnrow\textrm{Euler Characteristics } & -120 & -116 & -112 & -114&-108 \\
\hline
\end{tabular}$$
\end{thrm}
     
\subsubsection{\textbf{CY 3-folds with low codimension components}} In 4  cases  the Hilbert series of our codimension 4 CY 3-fold matches that of a known    codimension 3 Pfaffian CY 3-fold  but we find new deformation families in codimension 4. Among these,  three of them are orbifolds and one of them is a smooth CY 3-fold. In all orbifold cases, we have a weight   zero entry  in the syzygy matrix of the  Pfaffian CY 3-fold obtained after the
  projection. These Pfaffian CY 3-folds are degenerations of codimension 2 complete intersection CY 3-folds and the codimension 4 families are unprojections of these degenerations.
 
These 3 cases are  no. 10,  22, and 24 in Table \ref{main table} and have the same  Hilbert series of codimension 3 CY 3-folds with GRDB ID 745, 752, and 757 respectively on \cite{GRDB}.     We present a detailed analysis of no. 22 in Section \ref{CY3!codim2} and present the following theorem as a representative result.

\begin{thrm}\label{Th!low-codim} \label{Th!codim3}Let \(X\) be a Calabi--Yau 3-fold given with  the Hilbert series  
\begin{align}
P_{22}(t)=\frac{1- t^4- 2t^5-3t^6+3t^8 +4t^9+\cdots+t^{18}}{(1-t)^3(1-t^2)(1-t^3)^3(1-t^4)}=1+3t+ 7t^2+ 16t^3+\cdots
\end{align}
having  three  orbifold points of type \(\frac13(1,1,1)\) and \(D^3=2\). Then there exist four   deformation families realizing the Hilbert series \(P_{22}(t)\). Two of them are Tom, one is Jerry, and one is a codimension 3 Pfaffian family.  
 $$\begin{tabular}{|c|cccc|}\hline
\textrm{Family} & $\Tom$ & $\Tom$& $\Jer$& codim 3\\\hline
\evnrow\textrm{Euler Characteristics }e(X) & -110 &   -106 & -102&-94 \\
\hline
\end{tabular}$$
\end{thrm}         

\subsubsection{\textbf{Smooth  CY 3-folds in \(\w\bP^7\)}}
We find 4 families of smooth CY 3-folds among 23 working cases. The classification of smooth CY 3-folds in \(\PP^7\) has been conjecturally completed in \cite{CGKK}. The classification of all possible Betti tables for such CY 3-folds has been obtained in \cite{SSY}. Here we show the existence of three smooth CY 3-folds in strictly weighted projective space \(\w\bP^7\) and the other one which lives in straight \(\bP^7\), also appeared in \cite{CGKK}.   For  CY 3-folds corresponding to the  Hilbert series of  CY 3-fold  no. 5, we show the existence of a new
deformation family of the  codimension 3 Pfaffians CY 3-fold with ID 925 on \cite{GRDB}. We compute the Hodge numbers for both of them by using their explicit equations and they have different  Euler numbers:   they lie in two different components of the Hilbert scheme.  We prove the following theorem in Section \ref{CY3!smooth}.

\begin{thrm}\label{Th!smooth}The Hilbert series no. 1,3, 5 and 6 in Table \ref{main table} can be realized by a  smooth Calabi--Yau 3-fold in \(\pxp\) format. The no. 5 gives a new deformation family of a known   smooth Pfaffian  Calabi--Yau 3-fold of codimension 3.   
\end{thrm}    

\subsubsection{\textbf{Unprojections of non nodal Tom  formats  }}
In  the study  of \(\bQ\)-Fano 3-folds in \cite{BKQ},  each \(\pxp\)  family having a numerical type-I center admits a Gorenstein projection to a  nodal Pfaffian Fano 3-fold. For  CY 3-folds, this does not hold, despite the existence of  a numerical type-I center.  Further analysis of the Pfaffian 3-fold shows that none of the Tom and Jerry ansatz give a nodal degeneration.   In all such cases, we are only  able to find a single deformation family of Calabi--Yau 3-fold, in a $\pxp$ format.  There are 3 such cases; no. 14, 17, and 23 and we describe  the no. 14 in Section \ref{No!TJ}.  

\subsubsection{\textbf{Unprojection of a CY 3-fold with  curve of singularities}} 
All known cases of codimension 4 Fano and CY 3-folds with isolated orbifold points, having a  numerical type-I center, give a Gorenstein projection  to a  Pfaffian 3-fold  that at worst contains only isolated orbifold points. However, a CY 3-fold corresponding to the Hilbert series \(P_{20}\) admits a Gorenstein projection to Pfaffian CY 3-fold with a  \(\frac12(1,1)\) curve of singularities and we discuss this case in Section \ref{CY3!curve}.  

\subsection*{Acknowledgement}
We wish to thank Gavin Brown for many helpful discussions during this project. We are also grateful to an anonymous referee for their comments that helped us improve the exposition of the article.  SM would like to thank  Lahore University of Management Sciences for a doctoral studentship and  MIQ was supported by a research grant number INSS2308 of the Interdisciplinary Research Center of Intelligent and Secure Systems at the King Fahd University of Petroleum and Minerals, Dhahran, Saudi Arabia.

\section{Preliminaries}
\subsection{Notations and conventions}
\begin{itemize}
\item We work over the field of complex numbers \(\CC\).
\item A matrix of the form 

$$\begin{pmatrix}
m_{12} &m_{13} &m_{14}&m_{15}\\
&m_{23} &m_{24} &m_{25}\\
& &m_{34} &a_{34}\\
&&&m_{45}
\end{pmatrix}$$
represents the \(5\times 5\) skew-symmetric matrix or weights of such matrix. We omit the diagonal and lower triangular parts.
\item \(X, Y\) and \(Z\) represent varieties in codimension $4,3,$ and $2$ respectively.  
\item \(X_i\) for \(1\le i\le 5\) represents a Tom and \(X_{ij}\) represents
a Jerry codimension 4 CY 3-fold. 
\item A capital letter with subscript, for example \(F_i\), represents
a generic form of degree \(i\).

\end{itemize}

\subsection{Type-I unprojections and  Tom \& Jerry } 

We recall the notion of type-I unprojection, also known as Kustin--Miller
unprojection  \cite{KM,Ki}. Geometrically, it uses low codimension varieties
to construct varieties in higher codimension. 
\begin{dfn}Let $\bI({\Ybar})=(g_1,\ldots,g_s)$ be the ideal of a projectively Gorenstein variety  \(\Ybar\subset \PP(a_1,\ldots,a_n)\) of codimension $c$  and let $\bI(D)=(h_1,\ldots,h_t)$ be the ideal of a codimension $c+1$ variety $D$ contained in $\Ybar$. Then there exists a rational function $\nu$ on $\Ybar$
such that it has a pole along $D$, that is:
\begin{align*}
&\nu=\frac{l_1}{h_1}=\cdots=\frac{l_t}{h_t}; &\text{where }l_i\text{ are weighted
homogeneous polynomials.}
\end{align*}
The \textit{type-I unprojection} of   \(\Ybar\) is a projectively Gorenstein variety  $X$ of
codimension \(c+1\) defined by  $$\bV(g_1,\ldots,g_s,\nu h_1-l_1,\ldots,\nu h_t-l_t),$$in $\bP(a_1,\ldots,a_n,a_{n+1})$ where $a_{n+1}$ is the weight of $\nu$.
\end{dfn}
In codimension three,  most of the interesting classes of  Gorenstein varieties appear
as the vanishing locus of  the maximal Pfaffians of a $5\times 5$ skew-symmetric matrix \(M\). The type-I unprojection requires a specialization of  the equations of the Pfaffian variety to contain the divisor \(D\).     Tom and Jerry matrices are those specializations that provide the containment of the divisor
\(D\).

\begin{dfn}Let $J$ be the ideal of a  divisor \(D\)
of a codimension 3 Pfaffian variety.   Then for \(1\le i\le 5\), we define the \(i\)-\textit{th Tom matrix} $\Tom_i$  to be the matrix where the six entries $a_{jk}\in J$ for all $j,k\neq i$
and  the remaining four entries of the $i$th row and column are free choices.
For a pair \(ij: 1\le i,j\le 5\), we define the \(ij\)-\textit{th Jerry matrix}  $\Jer_{ij}$ to be the  matrix where the seven entries $a_{kl}\in J$ if either $k$ or $l$ equals $i$ or $j$. The rest of the three entries are free choices.\end{dfn}
 For example, if \(j_{ij}\) denote the entries in the ideal \(J\) and \(f_{ij}\) are free
entries then the  $\Tom_3$ and $\Jer_{45}$ matrices are given below.
\[\Tom_3=
\begin{pmatrix}
      j_{12}   &f_{13}  &j_{14}  &j_{15}\\
            &f_{23}  &j_{24}  &j_{25}\\
                       &&f_{34}  &f_{35}\\
                               &&    &j_{45}

 \end{pmatrix}\, \,  \, \,\ \ 
\Jer_{45}=
\begin{pmatrix}
         f_{12}   &f_{13}  &j_{14}  &j_{15}\\
                  &f_{23}  &j_{24}  &j_{25}\\
                       &   &j_{34}  &j_{35}\\
                               &&  &j_{45}

 \end{pmatrix}
\]


\subsubsection{\textbf{Computing number of nodes}}

To calculate the number of nodes on the  Tom or Jerry degeneration
$\Ybar$ of a codimension 3 Pfaffian variety $Y$ containing a  complete
intersection divisor \(D\), we can compute the number of
nodes either by computer algebra or by hand using the following formula of Brown--Kerber--Reid \cite{BKR}.
\begin{lema}\cite[Sec. 7]{BKR}\label{no. of nodes}
Consider the  skew-symmetric matrix  \(M\) 
\begin{align*}
M=\begin{pmatrix}
m_{12} &m_{13} &m_{14} &m_{15}\\
&m_{23} &m_{24} &m_{25}\\
&&m_{34} &m_{35}\\
&&&m_{45}
\end{pmatrix},
\end{align*}
corresponding to a codimension 3 Pfaffian variety \(Y\).
Let   \(\bI(D)=\langle x_1,...,x_4\rangle\) be an ideal of a divisor $D$ in \(Y\). 
If $M$ is in \emph{Tom}$_i$ format   then the four entries of the $ith$ row and column are free choices and provide a syzygy
\[
\Sigma_i=\sum_{\substack{k=1\\k \neq i}}^{5} m_{ik}\emph{Pf}_k.
\]
Let $d_k=\emph{deg }x_k$, $a_j=\emph{deg Pf}_j$ and $\sigma=\emph{deg }\Sigma_i$. Then the number of nodes on the divisor is given by the coefficient of the $h^2$ term in the expansion of
\[ \frac{\prod_{k=1}^4(1-d_kh)(1-\sigma h)}{\prod_{j\neq i}(1-a_jh)} .\]
If $M$ is a \emph{Jer}$_{ij}$ matrix for the ideal  $\bI(D)$, then the entries in the $ith$ and $jth$ rows and columns belong to the ideal. These give three syzygies $\Sigma_l$ of degree \(\si_l\) and another syzygy $T$ of degree $t=\emph{adjunction number}-\emph{wt }m_{ij}$. Then the number of nodes is again the coefficient of the $h^2$ term in the expansion of 
\[
\frac{\prod_{k=1}^4(1-d_kh)\prod_{l=1}^3(1-\sigma_l h)}{\prod_{j=1}^5(1-a_jh)(1-th)}
.\]

\end{lema}
 \subsection{Computing Hodge numbers of \( \PxP\) CY 3-folds}
\label{S!Hodge}
The most important invariants of a CY 3-fold \(X\) are its  Hodge numbers.  As all our varieties are quasismooth, so by \cite{steen} one can define a pure Hodge structure on them and    compute Hodge numbers  \(h^{p,q}\) as in the smooth case. For a CY 3-fold \(X\), computing \(h^{1,1}(X)\) and \(h^{2,1}(X)\) is sufficient to determine its Hodge diamond

\begin{equation}
\begin{array}{ccccccc} &&&1&&&\\ &&0&&0&& \\ &0&&h^{1,1}&&0& \\ 1&&h^{2,1}&&h^{2,1}&&1\\&0&&h^{1,1}&&0& \\ &&0&&0&&\\ &&&1&&&
\end{array}.\end{equation} 
 The Euler charactersistic \(e(X)\) is well known to be equal to \(2\left(h^{1,1}(X)-h^{2,1}(X)\right).\)

 We compute both  Hodge numbers of \( \PxP\) CY 3-folds \(X\) by using  \cite[Theorem 2.8]{DNFF} and implementation  the  \texttt{Versal Deformations} package of ~\cite{Ilten} in   the computer algebra system Macaulay2 \cite{M2}.  For each \(\PxP\) CY 3-fold we get  \(h^{1,1}(X)=2\); expected  due to \(\p2xp2\) being a Picard rank 2 variety. For other Tom and Jerry we can easily compute the equations by using the Papadakis method \cite{papadakis}. It is straightforward to  compute the Euler characteristics of other families but the explicit computation of \(h^{1,1}\) and \(h^{2,1}\)  by using computer algebra turns out to be computationally expensive as the equations in those cases are not as simple as those for \(\pxp\) families.   

All our families of CY 3-folds \(X\) have canonical orbifold points and they admit a crepant resolution to a smooth CY 3-fold \(\hat X\). The  resolutions of these orbifold points only change \(h^{1,1}(X)\), as described by Roan and Yau \cite[Sec. 2]{RY}. One can compute \(h^{1,1}(\hat X)\) by computing the  toric resolution of the singularities \(\frac 1r(a,b,c)\) \cite{YPG}.   We come across orbifold points with \(r\) equal to \(3,5\) and \(7\)  whose resolution adds \(1,2 \) and \(3\) to \(h^{1,1}(X)\) respectively by using Ron--Yau ~\cite{RY}.

\section{General Scheme of constructions}

\label{S!GenSch}
In this section, we describe a general scheme of the proofs for all cases appearing in Table \ref{T!Summary}.
\begin{enumerate}
\item[\textbf{Step 1:} ] \textbf{Existence of \(\pxp\) CY3.} We consider a candidate Calabi--Yau 3-fold in \(\pxp\) format, obtained using an algorithmic approach of \cite{QJSC,BKZ}. We prove the existence of the given canonical Calabi--Yau 3-fold \(X\) with the given invariants and singularities by checking the intersection of the orbifold strata of $\w\bP^7$  with $X$. We establish the quasismoothness either by explicitly studying base locus (see \cite{QMOC}) or by using the computer algebra system \textsc{Magma}\footnote{The code for each example can be accessed at https://github.com/QureshiMI/P2-x-P2-Quasismooth}. If the candidate \(X\) is  quasismooth then we
compute  Hodge numbers of \(X\) and its crepant resolution \(\hat X\), following
Section \ref{S!Hodge}.  Moreover,  If \(X\) contains a numerical type-I center, we move to step 2. 
\item[\textbf{Step 2:} ]\textbf{Special \(\pxp \) CY3.} Consider the embedding of $$X\subset \w\bP^7$$ 
of a Calabi-Yau 3-fold in codimension 4  given by  $\p2xp2$  format obtained at step 1. Let the weights of the variables \(x_{1},\ldots,x_7,z\) be given by  \(a_1,\ldots,a_7,r\) and let  \(P_Z:=\frac1r(a_5,a_6,a_7) \) be a type-I center, after necessary rearrangement of indices.  Then we deform the  equations of $X$ that can be expressed as the $2\times 2$ minors of a $3\times 3$ matrix:

\begin{equation}
\label{X!pxp}
\begin{pmatrix}
x_1 &x_2 & M\\
x_3 &x_4 & N\\
K &L &z
\end{pmatrix}
\end{equation}
where  $K,L,M,N$ are general forms. These are total 9 equations describing \(X\) in \(\bP(a_1,\ldots,a_7,r)\). In all cases,  we verify that  a presentation of the type
\eqref{X!pxp} gives a quasismooth \(\pxp\) CY 3-fold \(X\).

 \item[\textbf{Step 3:} ]\textbf{Projection from Type-I center.} In this step we perform a Gorenstein projection from the orbifold point  $P_{z}$, which is a type-I center. Geometrically it  means that the image of the point \(P_{z}\)  under projection is a weighted
projective plane $\bP(a_5,a_6,a_7)$ that is projectively normal in \(\Ybar \subset \bP(a_1,\ldots,a_7)\).  The projection is algebraically equivalent to excluding  4 out of 9 equations containing the variable  $z$ and the remaining 5 equations give a Calabi--Yau 3-fold $\Ybar\subset \bP(a_1,\ldots,a_7)$  in codimension 3 containing the divisor \(D=\bV(x_1,\ldots,x_4)\). The equations of $\Ybar$ can be described by  the $4\times4$
Pfaffians of the following skew-symmetric matrix,  having  zeros in the $m_{23}$ and $m_{45}$ positions 
\[
\begin{pmatrix}
K &L &M &N\\
&0 &x_1 &x_2\\
&&x_3 &x_4\\
&&&0
\end{pmatrix}.
\]
The homogeneity of the Pfaffian equations determines the degrees of the zero entries.
Without loss of generality we have written the equations of \(\Ybar\) in the  \(\Tom_1\) format, though we can write  in any of the   \(\Tom_i\) formats for \(i\in\{1,\ldots,5\}\). We check the singularities of \(\Ybar\)  and if \(\Ybar\)  contains only nodes that lie on \(D\) then we conclude that our \(\pxp\) CY 3-fold is of \(\Tom_1\) type.

 If \(X\) contains multiple type-I centers then we perform the projection from each type-I center point and continue with subsequent steps accordingly. There are two cases with two type-I centers that exhibit a similar phenomenon to \cite{BKR}, i.e. giving the same number of Tom and Jerry families with an equally  spaced number of nodes. 

\item[\textbf{Step 4:} ] \textbf{Deformation to Pfaffian CY 3-fold \(Y\).} In this step  we deform the entries of the \(\Tom_1\) matrix into a generic
Pfaffian matrix  to study the corresponding  codimension three CY 3-fold $Y\subset\bP(a_1,\ldots,a_7)$. This provides the resolution of nodes and
\(D\) is not a divisor of \(Y\). The equations of $Y$ are maximal Pfaffians of 
\[
M=\begin{pmatrix}
m_{12} &m_{13} &m_{14}&m_{15}\\
&m_{23} &m_{24} &m_{25}\\
& &m_{34} &m_{34}\\
&&&m_{45}
\end{pmatrix}
\]
where the  \(m_{ij}\) generic forms of the appropriate degrees in the weighted polynomial ring \(k[x_1,\ldots,x_7]\). For simplicity, we usually take the
entry in \(M\) to be \(x_i\) if \(\deg(m_{ij})=\wt(x_i)\).
If all the weights of \(M\) are positive,  except one example,  
  we land into one of the codimension three CY 3-folds constructed  in \cite{BKZ}, listed on \cite{GRDB}. In one case \(Y\)
turns out to be a CY 3-fold with a curve of singularities containing a dissident
singular point,  which is in itself a new model of CY 3-folds. If one of
the weights of \(M\) is non-positive then  it can not be quasismooth
by \cite[Proposition 2.7]{BKZ}.   
\item[\textbf{Step 5:} ] \textbf{Studying all Tom and Jerry degenerations of \(Y\).}  At this stage, we proceed with an additional 4 Tom and 10 Jerry degenerations of $Y$. We choose the entries $m_{ij}$ in such a way that the matrix takes the form of one of the Tom and Jerry matrix format having  $D$ lying in the corresponding degeneration $\overline{Y}$. For each    Tom and Jerry degeneration $\overline{Y}$ which contains only nodes as its singularities
and all of them lie on \(D\), there exists a corresponding  codimension four quasismooth CY 3-fold in  \(\bP(a_1,\ldots,a_7,r)\), with the same Hilbert series, orbifold points, etc. but distinct Euler characteristic due to the following lemma. \end{enumerate}
\begin{lema}\cite{clemens,Miles-CY}\label{L!EulNum} Let \(X\) and \(Y\) be CY 3-folds linked
via the diagram given below, then \(e(X)=e(Y)+2N-2\) where \(N\) is the number of
nodes on the divisor \(D.\)
\end{lema}
We  summarize our  generic scheme of study  in the following diagram:
\begin{center}

 \begin{alignat*}{3}
 \eqmathbox[L]{X\subset \text{w}\bP^7 } & & & \eqmathbox[R]{\hspace{15mm}X_{\text{T}_i/\text{J}_{ij}}\subset \text{w}\bP^7}\\
  \phantom{ = } \eqmathbox[L]{\dashdownarrow \scriptsize\text{T}_1} & & &\phantom{ + } \eqmathbox[R]{\hspace{15mm}\scriptsize\text{T}_i \dashdownarrow \scriptsize\text{J}_{ij} }\\
  D\subset \Ybar_{T_1} & \hspace{5mm}\rightarrow  &\hspace{5mm}Y &\hspace{5mm}\stackrqarrow{\text{T}_i}{\text{J}_{ij}}\eqmathbox[R]{\hspace{-2mm}D\subset
  \overline{Y}}
\end{alignat*}
\label{Figure 1}
\end{center}
In general the unprojection from a nodal Tom or Jerry family \(X \dashrightarrow  \Ybar\) factorizes through a  \(D\)-ample small resolution of nodes \(\Ybar \rightarrow \widetilde Y\), followed by a contraction \(\tilde Y \to X\) of the divisor \(D\).

\section{Calabi--Yau 3-fold with five distinct Tom and Jerry families}
\label{S!5TJ}
In the literature, all known examples of codimension 4 $\bQ$-Fano 3-folds contain at most 4 distinct  types of Tom and Jerry families, see \cite{BKRbigtable}. However, we show that the Hilbert series $P_7(t)$,  given in Table \ref{main table}, can be realized by five distinct Tom and Jerry families. The following provides the proof  of  Theorem \ref{Th!TJ5}.


\subsection{\(\pxp\) CY 3-fold family}We start with  a  weighted $\p2xp2$ variety $\cF$\[\w(\pxp)\into\bP(1^2,2^5,3^2)\] with the  weight matrix
\[
\begin{pmatrix}
1 &1 &2\\
2 &2 &3\\
2 &2 &3
\end{pmatrix}.
\] 
As \(\cF \) does not contain a singular locus of dimension 3, it is wellformed and  the  canonical divisor class is \(K_{\cF}=\Oh(-6)\).  We construct a triple projective cone over \(\cF\), by adding three new variables in the ambient ring which are not involved in the defining equations of $\w(\pxp) $,   giving a codimension four $\pxp$ \(7\)-fold \[\cC^3\cF \into\bP(1^5,2^5,3^2)\text{ with } K_{\cC^3\cF}=\Oh(-9).\] 
Then a  weighted complete intersection of \(\cC^3\cF\) with three general quadrics and one general cubic  gives a  3-fold 
\begin{align*}
X=\cC^3\cF\cap (2)^3\cap (3)\subset\bP(1^5,2^2,3):=\PP(x_1,\ldots,x_5,y_1,y_2,z),
\end{align*}
with trivial canonical class \(K_X=\Oh_X(-9+2.3+3)=\Oh_X\).  

The 3-fold \(X\) has the Hilbert series \(P_7(t)\) and  the Hilbert numerator for the ambient \(\pxp\) variety and \(X\) is the same, as it propagates the free resolution of their defining ideals. 
 The equations of $X$ are given by
\begin{equation}\label{eq!TJ5}
\bigwedge^2\begin{pmatrix}
x_1 &x_2 &A_2\\
y_1 &y_2 &B_3\\
C_2 &D_2 &z
\end{pmatrix}=0 
\end{equation}
where  $A_2, B_3, C_2, D_2$ are general forms.\\
Equivalently, we can describe \(X\) as a regular pullback from  \(\bP(1^5,2^2,3)\), i.e. the  entries of the matrix in \eqref{eq!TJ5} are weighted homogeneous forms in variables \(x_1,\ldots,x_5,y_1,y_2,z.\)
We can easily check that \(X\) does not intersect with the orbifold locus of weight 2. The orbifold point of weight 3 lies on $z$ and can be easily shown to be of type \(\frac13(1,1,1)\) by using the implicit function theorem.

\paragraph{\textbf{Quasismoothness}}
We prove the quasismoothness either by using a  version of Bertini's theorem
that if  a hypersurface $X \subset \bP(a_1,\ldots,a_n)$ is a general element
of a linear system $L = |\mathcal{O}(d)|$, then the singularities (non-quasismooth points) of
$X$ may only occur on the reduced part of the base locus of $L$  or by  using computer algebra. 
We explain one case in detail by using both approaches. Let \begin{align*}
V_1=\cC^3\cF&\subset\bP(1^5,2^2,3)
\end{align*}
be the triple  projective cone over $\cF$ with vertex $\bP^2$. Consider $V_2\subset V_1$, the locus of a general cubic. Then its base locus is given by
\begin{align*}
\text{Bs}(|\cO(3)|)&=\bV(\text{degree 3 weighted homogeneous monomials})\\
&=\bV(x_i^3,x_iy_j,z)=\bP(2^5).
\end{align*}
The linear system of the cubic has base locus $V_1\cap\bP(2^5)$. This is two-dimensional in $\bP(2^5)$ and   $V_2$ is quasismooth away from this locus.

Now let $X=\cC^3\cF\cap(2)^3\cap(3)\subset V_2$ be the general complete intersection of three quadrics. Then
\begin{align*}
\text{Bs}(|\cO(2)|)&=\bV(\text{degree 2 weighted homogeneous monomials})=\bV(x_i^2,y_j)
\end{align*}
which is a coordinate point of the variable \(z\) of degree 3 and \(X\) is already shown to be quasismooth on this point of type $\frac{1}{3}(1,1,1)$. Therefore, $X$ is a quasismooth CY 3-fold.

If the base locus is more complicated,  we can prove quasismoothness using computer algebra. We write down explicit equations over the rational
numbers of the given CY 3-fold and show that the affine cone of \(X\) is smooth by using  \textsc{Magma}~\cite{magma}. The challenge, in this case, is to come up with a sparse representation that can be handled easily by computer algebra.  For this example, the \(2 \times 2\) minors of the following give a quasismooth family.
$$\footnotesize{\begin{pmatrix}
x_1 &x_2 &y_2+y_1+x_4^2\\
y_1 &y_2 &z+x_1^3+x_2^3+x_3^3+x_4^3+x_5^3\\
x_5^2+x_4^2+x_1x_2+x_3x_2+x_1x_3+y_2 &x_3^2+x_2x_4+x_1x_5+x_{4}^2+x_5^2 &z \end{pmatrix}}$$
We compute Hodge numbers of \(X\) by using these equations in Macaulay2 \cite{M2}. We get  \(h^{1,1}(X)=2\) and \(h^{2,1}(X)=62\), giving Euler characteristics \(e(X)=-120 \). As \(X\) contains only one orbifold point of type \(\frac13(1,1,1)\), the crepant resolution \(\hat X\to X\) is a smooth CY 3-fold with \(h^{1,1}(\hat X)=3\) and \(h^{2,1}(\hat X)=62\). The Hodge pair \((3,62)\)  appears neither in the Kreuzer--Skarke \cite{KS} list nor in
the Green--H\"ubsch--L\"utken  \cite{CICY-List}, and is probably a  a new
Hodge pair. 
\subsection{Projection from Type I center}   We perform a Gorenstein projection from the type-I center point \(\frac13(1,1,1)\) to a Pfaffian CY 3-fold.  The projection from this point is a CY 3-fold \(\Ybar\) containing the divisor \(D=\bV(x_1,x_2,y_1,y_2)\) and it is given by the maximal Pfaffians of the  matrix  
\begin{align}\label{wMy}
\begin{pmatrix}
C_2 &D_2 &A_2 &B_3\\
&0 &x_1 &y_1\\
&&x_2 &y_2\\
&&&0
\end{pmatrix}.
\end{align}
The matrix is given in \(\Tom_1\) format, i.e. the entries other than the first row or column are in \(\bI(D)\): the ideal of the divisor \(D\).  Then by using  Lemma \ref{no. of nodes} or computer algebra, we can show that \(\Ybar\) contains 10 nodes and all of them lie on \(D\). The unprojection of the pair \((\Ybar, D)\) gives the \(\pxp\) CY 3-fold described by the matrix \eqref{eq!TJ5}.  

\subsection{\textbf{ Pfaffian CY 3-fold \(\Y\).}} The homogeneity of the Pfaffian equations gives deg $m_{23}=1$ and deg $m_{45}=2$. A deformation \(\Ybar \to Y\ \) provides a resolution of nodes and gives a quasismooth Pfaffian CY 3-fold 
$$Y_{3^2,4^3}\subset\bP(1^5,2^2)$$
 with the weight matrix 
\[
\begin{pmatrix}
2 &2 &2 &3\\
&1 &1 &2\\
&&1 &2\\
&&&2
\end{pmatrix}.
\]
In fact, this is the family of  CY 3-fold with GRDB ID. 928 on \cite{GRDB}, constructed in \cite{BKZ}. Then by using Lemma \ref{L!EulNum}
$$e(Y)=-120-2\times10+2=-138.$$

\subsection{Further  Tom and Jerry families} 
In this section, we find other deformation families of \(X\) by studying other  possible 4  Tom  and 10  Jerry degenerations \(Y \rightsquigarrow \Ybar_{\Tom_i,\Jer_{ij}}\).
 If the corresponding degeneration is at worst nodal, then its unprojection will give a  quasismooth CY 3-fold in codimension 4 with Hilbert series \(P_7(t)\) and a single orbifold point \(\frac13(1,1,1)\). 
\subsection*{Tom CY 3-fold $X_{2}$}
To construct $\Tom_2$ we need the entries of the 2nd row and column of the weight matrix to be the general forms and all other entries will be in  $\bI(D)$.
For example,\begin{align*}
\begin{pmatrix}
K_2 &y_1 &y_2 & \langle x_4, x_5\rangle_3\\
&L_1 &M_1 &N_2\\
&&x_1 &\langle x_2, y_2\rangle_2\\
&&&\langle x_2, y_1\rangle_2
\end{pmatrix},
\end{align*}
where $K_2, L_1, M_1$, and $N_2$ are the general forms of the  degrees equal to their subscripts and the rest of the entries are  taken to be in $\bI(D)$, after performing row and column operations on the more general entries in \(\bI(D)\). Then  the maximal Pfaffians of the matrix give a nodal CY 3-fold with  12 nodes on $D$, which we compute by using Lemma \ref{no. of nodes}.
We have $$Y_{3^2,4^3}\subset\bP(1^5,2^2), D=\bP(1,1,1),$$ then $d_i=1,1,2,2$, $a_j=\text{wt Pf}_j$ and $\sigma_1=\text{wt }\Sigma_1$ where $$\Sigma_1=C_2\text{Pf}_2+D_2\text{Pf}_3+A_2\text{Pf}_4+B_3\text{Pf}_5.$$
Then 
\begin{align*}
f(h)&=\frac{\prod_{i=1}^4(1-d_ih)(1-\sigma_1h)}{\prod_{j\neq1}(1-a_jh)}=\frac{(1-h)^2(1-2h)^2(1-6h)}{(1-3h)(1-4h)^3}.
\end{align*} 
The number of nodes is then the coefficient of the $h^2$ term in the expansion of $f(h)$, which equals $10$. The Euler number  \(e(X_{{2}})=-138+2\times12 -2=-116\).

Notice that due to the symmetry in the weight matrix, interchanging the 2nd row with the 3rd or 4th row gives back the same weight matrix. Therefore, $\Tom_2$ is equivalent to $\Tom_3$ and $\Tom_4$, and all of them give 12 nodes on $D$.
Furthermore,  the 1st and 5th rows and columns also have the same weights so $\Tom_5$ is equivalent to $\Tom_1$ and has 10 nodes on the divisor \(D\).

\subsection*{Jerry CY 3-folds $X_{{12}},X_{{15}},X_{{23}}$} 
We find 3 more deformation families of \(X\) by using Jerry degenerations, and one of them is presented below.   To construct  \(\Jer_{12}\), we take entries \(m_{jk}\) with \(j,k =3,4,5 \) to be general forms and the entries where \( j\) or \(k\) equal to either 1 or 2 to be in the ideal of \(D\). For example, the following matrix is a presentation of \(Y\) in \(\Jer_{12}\) format, where the subscripts of \(K, L, M,\langle y_2\rangle\) and \( \langle x_1,x_2\rangle\) denotes their degree.     
\[
\begin{pmatrix}
y_1 &y_2 &\langle x_1, x_2\rangle_2 &\langle x_1, x_2\rangle_3\\
&x_1 &x_2 &\langle y_2\rangle_2\\
&&K_1 &L_2\\
&&&M_2
\end{pmatrix}
\]
 Then  $\Jer_{12}$ contains  14 nodes lying on the divisor \(D\), which we calculate from  \cite[Sec. 7]{BKR} given in Lemma \ref{no. of nodes}.
Therefore the  Euler characteristic is given by   \[e(X_{12})=-138+2\times 14-2=-112.  \]  Since the 2nd, 3rd and 4th rows and columns have identical weights we see that $\Jer_{13}$ and $\Jer_{14}$ are identical to $\Jer_{12}$ and also have 14 nodes on $D$.

For the case of \(\Jer_{15}\), we can take the following matrix 
\[
\begin{pmatrix}
y_1 &y_2 &\langle x_1,x_2\rangle_2 &\langle x_1,x_2\rangle_3\\
&K_1 &L_1 &\langle x_1,x_2, y_2\rangle_2\\
&&M_1 &\langle x_1,x_2,y_1\rangle_2\\
&&&\langle x_1,x_2,y_1,y_2\rangle_2
\end{pmatrix},
\]
and show that it contains  13 nodes on the divisor \(D\). Therefore its Euler characteristic is  \(e(X_{{15}})=-138+2\times13-2=-114\). 

The case \(\Jer_{23}\) can be represented by the Pfaffians of the matrix 
\[
\begin{pmatrix}
y_1 &y_2 &K_2 &L_3\\
&x_1 &x_2 &\langle y_2\rangle_2\\
&& 0 &\langle x_2, y_1\rangle_2\\
&&&M_2
\end{pmatrix}
.\]
It contains 16 nodes, so we get  \[e(X_{23})=-138+2\times16-2=-108. \]Following the same reasoning as before, we see that $\Jer_{24}$ and $\Jer_{34}$ are also equivalent to $\Jer_{23}$ because the 2nd, 3rd, and 4th rows and columns have the same weights, so each case has 16 nodes on $D$.

Lastly, we can see that $\Jer_{25}$ is equivalent to $\Jer_{35}$ and $\Jer_{45}$, and our calculations show that each of these has 14 nodes on the divisor. 

\section{CY 3-folds as unprojection of degenerations of complete intersections}
\label{CY3!codim2}

In this section, we discuss CY 3-folds in codimension 4 that are obtained as unprojection of
Pfaffians that are degenerations of codimension 2 complete intersections. They   exhibit  new deformation families in the Hilbert scheme of the  existing codimension 3 Pfaffian CY 3-folds
of Brown--Kasprzyk--Zhou \cite{BKZ}.  In the orbifold case, the Hilbert series   no.  10, 22, and 24 in Table \ref{main table}  matches with  the Hilbert series of CY
3-folds with  GRDB ID 745, 752, and 757 respectively. For  745 and 757, we
get 2 new deformation families and  for 752 we get 3 new deformation families. In the smooth case, the Hilbert series \(P_{5}\) matches that of the CY 3-fold with
GRDB ID 925. We discuss the deformation families of  CY 3-folds corresponding to Hilbert series \(P_{22}(t) \) that provides the proof of  Theorem \ref{Th!low-codim}.

\subsection{\(\pxp\) CY 3-fold}
Consider a weighted \(\pxp\) variety \(\cF\) with the weight matrix 
\[
\begin{pmatrix}
1 &2 &3\\
2 &3 &4\\
3 &4 &5
\end{pmatrix}.\]
It is well formed and its canonical divisor class is $K_{\cF}=\cO(-9)$. Taking a double projective cone over $\cF$ and taking a weighted  complete intersection with a quadric, a quartic and a quintic gives a 3-fold
\[
X=\cC^2\cF\cap (2)\cap(4)\cap (5)\subset\bP(1^3,2,3^3,4):=\PP(x_1,x_2,x_3,y,z_1,z_2,z_3,u)
\]
with trivial canonical class $K_X=\cO_X$. 
The Hilbert series of $X$ is 
\begin{align}
    \label{HS!C4}
P_{22}(t)=\frac{1- t^4- 2t^5-3t^6+3t^8 +4t^9+\cdots+t^{18}}{(1-t)^3(1-t^2)(1-t^3)^3(1-t^4)}=1+3t+ 7t^2+ 16t^3+\cdots .
\end{align}
The defining equations of $X$ are given by 
$$\bigwedge^2\begin{pmatrix}
x_1 &A_2 &z_1\\
y &z_2 &B_4\\
z_3 &u &C_5
\end{pmatrix}=0,$$
where  $A_2, B_4, C_5$ are general forms.

It can be shown that the singular locus of $X$ consists of three orbifold points, each of type $\frac{1}{3}(1,1,1)$ and $X$ is a quasismooth CY 3-fold, following  the last section.  As discussed in Section \ref{S!Hodge} we get   \(h^{2,1}(X)=
57\) and therefore \(e(X)=-110\).

 \subsection{Projection from type-I center}
A Gorenstein projection from the type-I center \(P_{z_1}\)
gives   a Pfaffian CY 3-fold \(\Ybar\) that is the vanishing locus of the
  $4\times 4$ Pfaffians of
\begin{align}\label{wMy1}
\begin{pmatrix}
    x_1 &A_2 &B_4 &C_5\\
    &0 &y &z_2 \\
    &&z_3 &u\\
    &&&0
\end{pmatrix},
\end{align}
containing the divisor \(D=\PP^2=\bV(y,z_2,z_3,u)\). The homogeneity of the Pfaffian equations forces the  deg $m_{23}=0$ and deg $m_{45}=6$.  Then by using Lemma \ref{no. of nodes} or computer algebra, we can show that \(\Ybar\) contains 22 nodes and all of them lie on \(D\), so its unprojection is a \(\Tom_1\) CY 3-fold \(X\) in  \( \pxp\)  format.

A  deformation   \(\Ybar\rightarrow\ Y\) of entries of $M$ gives a CY 3-fold  given by $$Y_{4,5,6,7,8}\subset\bP(1^3,2,3^2,4)$$ with weight matrix 
 
\begin{align}\label{WM-degen}
\begin{pmatrix}
1 &2 &4 &5\\
&0 &2 &3\\
&&3 &4\\
&&&6
\end{pmatrix}.
\end{align}
Notice that the matrix of the weights of $Y$ contains a weight zero entry, so by \cite[Prop. 2.7]{BKZ}, $Y$ cannot be quasismooth.
The degree 0 entry is a zero polynomial in the matrix \eqref{wMy1} but 
in matrix \eqref{WM-degen} it can be any constant.  One can describe $Y$
by the Pfaffians of\begin{align*}
\begin{pmatrix}
x_1 &y &u &K_5\\
&\alpha &L_2 &z_1\\
&&z_2 &M_4\\
&&&N_6
\end{pmatrix}, 
\end{align*}
where $K_5, L_2, M_4, N_6$  are general forms and $\al$ is a constant. The following are the Pfaffian equations of $Y\subset \bP(1^3,2,3^2,4)$:
\begin{align*}
\Pf_1&=\alpha N_6-L_2M_4+z_1z_2\\
\Pf_2&=yN_6-uM_4+K_5z_2\\
\Pf_3&=x_1N_6-uz_1+K_5L_2\\
\Pf_4&=x_1M_4-yz_1+\alpha K_5\\
\Pf_5&=x_1z_2-yL_2+\alpha u.
\end{align*}
 By varying  the entry \(m_{23}\) in the Pfaffian matrix to be equal to 1,
 $\Pf_5$ gives us $u=yL_2-x_1z_2$. This provides two syzygies that eliminate two of the four remaining Pfaffian equations. Taking this value of $u$, we get
\begin{align*}
\Pf_2&=yN_6-(yL_2-x_1z_2)M_4+K_5z_2\\
\Pf_3&=x_1N_6-(yL_2-x_1z_2)z_1+K_5L_2.
\end{align*}
Notice that $\Pf_2=y\Pf_1+z_2\Pf_4$ and $\Pf_3=x_1\Pf_1+L_2\Pf_4$. This means that we are left with only 2 equations of degree 5 and 6 so that $Y_{4,5,6,7,8}\subset \bP(1^3,2,3^2,4)$ is just a degeneration of $Y_{5,6}\subset \bP(1^3,2,3^2)$.
 Therefore, the degeneration \(Y_\alpha\to
\AA^1\) has generic fibre \(Z_{5,6}\) and the central fibre \(Y_0\) is a codimension
3 Pfaffian CY 3-fold. 

 The Euler characteristic \(e(Z_{5,6})=-152\)
is well known. Given that the projection from a \(\pxp\) CY 3-fold \(X\) has 22 nodes
and \(e(X)=-110\), one can see that \[e(X)=e(Z_{5,6})+2\times 22-2,\]in accordance
with Lemma \ref{L!EulNum}. So we conclude that \(X\) can be obtained from the unprojection
of the  \(\Tom_1\) format of the Pfaffian degeneration of \(Z_{5,6}\).

\subsection{Further  Tom and Jerry families} Now we study all possible Tom and Jerry degenerations to find the ones that are at worst nodal. This gives further two deformation families whose unprojection will give a quasismooth CY 3-fold in codimension 4 with Hilbert series $P_{22}(t)$ and three orbifold points  of type $\frac{1}{3}(1,1,1)$.
\subsection*{A Tom CY 3-fold $X_2$} To construct $X_2$ we consider the $\Tom_2$ format of the weight matrix \eqref{WM-degen}, i.e.
\[
\begin{pmatrix}
    K_1 &y &u &\langle z_2,z_3\rangle_5\\
    &\alpha &L_2 &M_3\\
    &&\langle z_2,z_3\rangle_3 & \langle z_2,z_3, u\rangle_4\\
    &&&\langle y,z_2,z_3\rangle_6
\end{pmatrix}
.\]
This  has 24 nodes on the divisor $D=\bP^2$ and its unprojection to codimension 4 is a CY 3-fold $X_2$ with Euler characteristic:
\[e(X_2)=-152+2\times24-2=-106.\]
\subsection*{A Jerry CY 3-fold \(X_{35}\)} If we take the \(\Jer_{35}\) format of \eqref{WM-degen}, then it can be described by the maximal Pfaffians of   \[
\begin{pmatrix}
K_1 &y &L_4 &\langle z_2,z_3\rangle_5\\
&0 &M_2 &\langle z_2,z_3\rangle_3\\
&&z_2 &u\\
&&&\langle y,z_3\rangle_6
\end{pmatrix}.
\]
We calculate that it contains 26 nodes lying on the divisor \(D=\PP^2\) and
its unprojection is a codimension 4 CY 3-fold with Euler characteristics:
\[e(X_{35})=-152+2\times 26-2=-102.\] 
\subsection{Codimension 3 family ID 752 }
Consider the CY 3-fold \(Y\) ID 752 on the the graded ring data base \cite{GRDB}, then its Hilbert series \[P_{752}(t)=\dfrac{1 - 2t^5 - 3t^6 + 3t^8 + 2t^9 - t^{14}}{(1-t)^3(1-t^2)(1-t^3)^3}=1 + 3t + 7t^2 + 16t^3 +\cdots ,\]
matches the Hilbert series \eqref{HS!C4} of codimension 4 deformation families. We show that this lies in a different deformation family than the
above three families. The weight matrix of \(Y\) is   \[
\begin{pmatrix}
2 &2 &3 &3\\
&2 &3 &3\\
&&3 &3 \\
&&&4
\end{pmatrix}, 
\]
  and we have the embedding \[Y_{5^2,6^3}\subset \PP(1^3,2,3^3)\]
  containing three orbifold points of type \(\frac{1}{3}(1,1,1)\), each of which are
type-I centers. The projection
from any one of these points is a complete intersection CY 3-fold  \[(\overline Z_{5,6}\supset D=\PP^2)\subset \PP(1^3,2,3^2).\] It contains 30 nodes that lie on the
divisor \(D\). The small resolution \(\overline Z \to  Z\) of nodes gives a  CY 3-fold \(Z_{5,6}\subset \PP(1^3,2,3^2)\) with \(e(Z_{5,6})=-152\).
Therefore by Lemma \ref{L!EulNum}, we have \(e(Y)=-152+2\times30 -2=-94\)
and \(Y\) is a different deformation family than the three codimension 4 Tom and Jerry families discussed above.
   \section{Smooth Calabi--Yau 3-folds}
\label{CY3!smooth}
Among our 23 working cases, we found four of these, no. 2, 3, 5, and 6 in Table \ref{main table}, to be smooth CY 3-folds, i.e. they are Calabi--Yau manifolds. They can also be obtained as  linear sections of the smooth  Fano 4-folds constructed  in \cite{QBAM}, but here  they are of independent interest.  The CY 3-fold no. 2 was given in \cite{CGKK} but other cases present new models of smooth  CY 3-folds in \(\w\bP^7\). Their degree \(D^3\) and Hodge numbers are listed in Table \ref{T!Summary}. In this section, we discuss no. 5 which also provides a proof of Theorem \ref{Th!smooth} by showing  that there are two deformation families of CY 3-folds corresponding to the Hilbert series \(P_5(t)\). 
\subsection{\(\pxp\) family}
Consider a weighted \(\pxp\) variety \[\cF\into\bP(1^4,2^4,3),\] with the weight matrix \[
\begin{pmatrix}
1 &1 &2 \\
1 &1 &2\\
2&2&3 
\end{pmatrix}.
\] Then the canonical class is \(K_\cF=\Oh(-5)\) as \(\cF\) is wellformed, i.e. does not contain a singular divisor. Indeed, \( \cF\) contains singularities coming from a weight 3 point and  weight 2 locus. We take two projective cones over \(\cF\) to get \[V:=\cC^2\cF\into \PP(1^6, 2^4,3) \text{ with } K_{\cC^2\cF}=\Oh(-7).\] Then the complete intersection of \(V\) with a general cubic \(A_3\) and two quadrics \(B_2,C_2\) is a CY 3-fold \[X=V\cap A_3\cap B_2\cap C_2\into\PP(1^6,2^2):=\PP(x_1,\ldots,x_6,y_1,y_2).\] The equations of \(X\) are given by 
\[\bigwedge^2\begin{pmatrix}
x_1 &x_2 &y_1 \\
x_3 &x_4 &y_2\\
B_2&C_2&A_3\end{pmatrix} \] where $A_3,B_2,C_2$ are general forms.

The weight 2 locus is \[\VV(B_2y_1,B_2y_2,C_2y_1,C_2y_2)\into \PP(y_1,y_2),\] which is an empty set, so \(X\) is a smooth CY 3-fold. The Hodge numbers are \(h^{1,1}=2\) and \(h^{2,1}=50\), giving the Euler characteristic to be \(e(X)=2(2-50)=-96.\)   
\subsection{Pfaffian family}  The Hilbert series of the CY 3-fold with GRDB ID. 925 matches that of \(P_5(t)\). The ID. 925 is a smooth Pfaffian CY 3-fold \(Y\) that can be described in \(\PP(1^6,2)\)  by the maximal Pfaffians of the skew-symmetric matrix 
\[\begin{pmatrix}
x_1 &x_2 &x_3 &x_4\\
&y &y+x_2^2+x_5x_6+x_5^2 &y+x_1^2+x_3x_5+x_6^2\\
&&x_6x_4+x_5^2-x_1x_3+x_4^2 &x_5^2-x_6^2+x_3x_4+y\\
&&&x_5^2+x_1x_4
\end{pmatrix}
\] where $\wt \ x_i=1$ and $\wt\ y=2$.
 We compute the Hodge numbers using this  sparse presentation of  equations of \(Y\), which gives \(h^{1,1}(X)=1\) and \(h^{2,1}(X)=59\), with Euler characteristic given by  \[e(Y)=2(1-59)=-116. \]This shows that the two families are topologically distinct.  $\square$

\section{CY 3-fold from non nodal  Tom unprojection}
\label{No!TJ}
In the case of $\QQ$-Fano 3-folds, it was noticed in \cite{BKQ} that a Gorenstein projection from a  numerical type-I center of every quasismooth $\pxp$ Fano 3-fold is a nodal Pfaffian Fano 3-fold $\Ybar$ that admits a resolution of nodes to a quasismooth Pfaffian $\QQ$-Fano 3-fold $\Ygen$. For the case of  $\p2xp2$ CY 3-folds, this no more holds. In the families no. 14, 17, and 23 in Table \ref{main table}) we discover that neither the image of projection $\Ybar$ is nodal nor the general Pfaffian CY 3-fold \(\Ygen\) is quasismooth. Moreover, all other Tom and Jerry degenerations of \(\Ygen\)  also contain
worse singularities than nodes on \(D\).  We briefly discuss the family no. 14 for the description.

Following  the ideas  from earlier sections, we can show that  $$X\subset\bP(1^4,2^2,3,5):=\PP(x_1,\ldots,x_4,y_1,y_2,z,u)$$ with a weight matrix 
\[
\begin{pmatrix}
    1 &1 &2\\
    2 &2 &3\\
    4 &4 &5
\end{pmatrix}  
\]
is a quasismooth CY 3-fold with two orbifold points: \(\frac13(2,2,2) \) and \(\frac15(1,1,3)\). The equations of \(X\) can be given  by 
\[\bigwedge^2
\begin{pmatrix}
x_1 &x_2 &A_2\\
y_1 &y_2 &z\\
B_4 &C_4 &u
\end{pmatrix}\ \ \ \ \ \ \text{where } A_2, B_4, C_4\text{ are general forms}.\]
The orbifold point \(\frac13(2,2,2)\) is not a type-I center, so we perform a projection from  a type-I center \(P_u:=\frac 15(1,1,3)\). The image of the projection is given by  the $4\times4$ Pfaffians of the following $5\times5$ skew-symmetric matrix:
\[
\begin{pmatrix}
    A_2 &z &B_4 &C_4\\
    &0 &x_1 &x_2\\
    &&y_1 &y_2\\
    &&&0
\end{pmatrix}
.\]
By the homogeneity of the Pfaffian equations, the degree of $m_{23}$ is zero, and the degree of $m_{45}$ is equal to 3. This a $\Tom_1$ presentation of the projection so   the entries in the first row and column are in the ideal of the divisor  $$D = \VV(x_1, x_2, y_1,y_2) = \PP(1,1,3).$$ However, notice that the pure power of the coordinate variable $z$ can appear as a multiple of a linear variable in only two  equations of  \(Y\), so $P_z$ is a point of higher embedding dimension in $Y$. Hence, $\Tom_1$ fails and \(\Ybar\) contains singularities that are worse than nodes on the divisor \(D\).

The deformation to \(\Ygen\) is also not quasismooth by \cite{BKZ} as the syzygy weight matrix contains a  zero entry. All other  Tom and Jerry degenerations of \(Y\)  also fail due to having points of higher embedding dimension. 
The CY 3-folds with Hilbert series No. 17 and 23 also propagate the same phenomena.   
\section{CY 3-folds  as unprojections of CY 3-folds with curve singularities}
\label{CY3!curve}

In all known cases of codimension four  Fano 3-folds with numerical type-I centers and   isolated orbifold points,  the image   Pfaffian Fano 3-fold also contains isolated orbifold points at worst. However, a  CY 3-fold No. 20 in Table \ref{main table} that contains an isolated point of type \(\frac17(1,2,4)\), admits a Gorenstein projection to a Pfaffian CY 3-fold \(Y\) with a curve of singularities \(C\), containing a  non-isolated orbifold point lying   on \(C\).  

Consider a  weighted \(\pxp\) variety $\cF\subset\bP(1^4,4^4,7)$ with weight matrix 
\[
\begin{pmatrix}
1 &1 &4\\
1 &1 &4\\
4 &4 &7
\end{pmatrix}. 
\]Suppose that $V:=\cC^2_{1,2}\cF$ is a double cone over $\cF$, where the cone variables are of degree 1 and 2.
Then a weighted complete intersection of \(V\) with 3 generic quartics \(X=V\cap \left\{ Q_i\right\}_{i=1}^3\) is a CY 3-fold  \begin{align*}
&X\subset\bP(1^5,2,4,7):=\bP(x_1,x_2,x_3,x_4,x_5,y,z,u).
\end{align*} 
The equations of $X$ can be described by the  $2\times 2$ minors of the below $3\times 3$ matrix 
\begin{align*}\label{eq!diss}
&\begin{pmatrix}
x_1 &x_2 &z\\
x_3 &x_4 &A_4\\
B_4 &C_4 &u
\end{pmatrix} &&\text{where } A_4, B_4, C_4\text{ are general forms. }
\end{align*}
One can show that \(X\) is a quasismooth CY 3-fold with a single orbifold point of type \(P_u:=\frac17(1,2,4)\) which is a type-I center.    
The Hodge number \(h^{2,1}= 89\), giving Euler characteristics to be \(e(X)=-174\). The projection from  $P_u$ is a  CY 3-fold \(\Ybar\) containing the divisor \[D=\PP(1,2,4)=\bV(x_1,x_2,x_3,x_4)\] that can be described in \(\Tom_1\) format by the maximal Pfaffians of 
\[
\begin{pmatrix}
  B_4 &C_4 & A_4 & z\\
  & 0 & x_3 & x_1\\
  && x_4 & x_2\\
  &&& 0
\end{pmatrix}.
\]
The homogeneity of the Pfaffian equations forces the degree of $m_{23}$ and $m_{45}$ to be equal to 1.   Then singularities  on $\Ybar$ turn out to be 4 nodes that all lie on  $D$. In this case the divisor \(D=\PP(1,2,4)\) is non-normal and contains a curve \(C\) of singularities of type \(\frac12(1,1)\)  and a dissident singular point \(\frac14(1,1,2)\) lying on \(C\). We get \(X\)  as an unprojection of a non-normal divisor \(D\) containing a curve of singularities and unprojection providing a resolution of singularities and  contraction of the divisor \(D\) to the point of type \(\frac17(1,2,4)\).  
     
  Then  a  deformation \(\Ybar_{\Tom_1}\rightarrow Y\) is a CY 3-fold  $$Y_{2,5^4}\subset\bP(1^5,2,4)$$ with a syzygy matrix \[
\begin{pmatrix}
    4 &4 &4 &4\\
    &1 &1 &1\\
    &&1 &1\\
    &&&1
\end{pmatrix}.
\]
The  orbifold locus of weight 4 is a single coordinate point \(P_z\).
If we fix \(z\) at position \(m_{15}\) and using weight one variables at position \(m_{23}, m_{24}\) and \(m_{34}\), we can show that it is  a nonisolated orbifold point of type \(\frac14(1,1,2)\), called a dissident singular point. The intersection of weight \(2\) with \(Y\) is one dimensional and  we can show that it is an orbifold curve of singularities of type \(\frac12(1,1)\) that contains the dissident point \(P_z\). By using Lemma  \ref{L!EulNum},
we get \(e(Y)=-174-8+2=-180.\)  

Further  analysis reveals the existence of one more deformation family of \(X_{23}\) as an unprojection of a nodal \(\Jer_{23}\) containing 6 nodes on \(D\) with \(e(X_{23})=-170\).
In total this gives rise to 2 distinct codimension 4 deformation families.

\appendix
\section{Table}
The following table encompasses the details of the  calculations performed on all the cases. The first 2 columns detail the embedding of a codimension 4 CY 3-fold, and the weight matrix of the corresponding \(\pxp\) family, and the orbifold points respectively. If \(X\) has a numerical type-I center then the column ``T-I.'' lists the  type of that point. The next column gives the equation degrees and the embedding of the corresponding Pfaffian CY 3-fold \(Y\) in \(\PP(a_1,\ldots, a_7)\). If \(Y  \) is in the GRDB \cite{GRDB} then it also lists its GRDB ID in the same column. The  column ``Syzygy WM'' lists the  weights of the Pfaffian CY\ 3-fold \(Y\). The last  column lists the working Tom and Jerry cases, alongside their number of nodes. The Hilbert numerator propagating the  codimension 4 graded free resolution  and the first ten coefficients \(p_m(t)\) of the Hilbert series $\sum_{m\ge 0}p_m(t)t^m$ are given below each corresponding entry of the table.
\begin{rmk} The weight matrices in the table  are rearranged to list the weights of the syzygy matrix in  increasing order in each row. So   the Tom family corresponding to the \(\pxp\) CY 3-fold is not always  \(\Tom_1\) as described in the generic scheme of analysis in Section \ref{S!GenSch} and followed in subsequent sections. The Tom family listed in bold corresponds to the \(\pxp\) CY 3-fold.   
\end{rmk} 
\begin{rmk} Families of  CY 3-folds   no. 10, 22, and 24 appear as codimension 4 deformation components of the Hilbert scheme  of already existing codimension 3 Pfaffian CY 3-folds. We list the GRDB ID of the corresponding Hilbert series of the codimension 3 CY 3-fold   \(Y\) that appeared in \cite{GRDB}  and \(N\) in ``TJ\#nd'' column denotes the nodes on the projection from \(Y\) to a
codimension 2 complete intersection CY 3-fold \(Z\).   
\end{rmk}

    \pagestyle{plain}

\begin{center}

\begin{longtable}{p{0.4cm}p{3.1cm}p{2cm}p{0.7cm}p{3cm}p{2.5cm}p{1.4cm}}

\hline

 No. & $X\subset \text{w}\bP^7$ & $\bP^2\times\bP^2$ WM &  T-I & \textsc{Grdb} ID \&\newline$Y\subset\text{w}\bP^6$ & Syzygy WM  & TJ \#nd \\

\hline

\endfirsthead

\hline

No. & $X\subset \text{w}\bP^7$ & $\pxp$ WM  &  T-I & \textsc{Grdb} ID/\newline$Y\subset\text{w}\bP^6$ & Syzygy WM  & TJ \#nd\\

\hline
\endhead

\hline
\multicolumn{7}{r}{\scriptsize \emph{Continued on next page}} \\
\endfoot

\hline
\endlastfoot

\oddrow
  {1}    & $\begin{array}{@{}l@{}} X_{2^3,3^6}\\\quad\subset\bP^7\end{array}$ & \footnotesize{$\begin{pmatrix}
1 &1 &1\\
1 &1 &1\\
2 &2 &2
\end{pmatrix}$}&&&& \\
 &\multicolumn{6}{l}{\footnotesize{Hilbert Numerator: \(1- 3t^2- 4t^3+ 12t^4- 4t^5 - 3t^6+t^8       \)}}\\
&\multicolumn{6}{l}{\footnotesize{Hilbert Coefficients: $ [ 1, 8, 33, 92, 202, 380, 643, 1008, 1492, 2112, \ldots ]$}}\\

 \evnrow {2}  &$ \begin{array}{@{}l@{}} X_{2^3,4^6}\\\quad\subset\bP(1^7,3)\end{array}$ & $\footnotesize{\begin{pmatrix}
1 &1 &1\\
1 &1 &1\\
3 &3 &3
\end{pmatrix}}$ &  1/3 & $Y_{2^3,4^2}\subset \bP^6$ & \footnotesize{$\begin{pmatrix}
-1 &1 &1 &1\\
        &1 &1 &1\\
                &&3 &3\\
                        &&&3
\end{pmatrix}$} & 
$\begin{matrix}
    \textbf{T$_3$}& 10\\
    J_{13}&13
\end{matrix}$\\  
\evnrow&\multicolumn{6}{l}{\footnotesize{Hilbert Numerator: \( 1- 3t^2+2t^3-6t^4+12t^5-6t^6+2t^7-3t^8+t^{10} \)}}\\
 \evnrow&\multicolumn{6}{l}{\footnotesize{Hilbert Coefficients: $ [1, 7, 25, 66, 141, 261, 438, 683, 1007, 1422,\ldots ]$}}\\


3 &  $ \begin{array}{@{}l@{}} X_{3^6,4^3}\\\quad\subset\bP(1^6,2^2)\end{array}$ & $\footnotesize{\begin{pmatrix}
1 &1 &1\\
2 &2 &2\\
2 &2 &2
\end{pmatrix}}$ \\
&\multicolumn{6}{l}{\footnotesize{Hilbert Numerator: \(1- 6t^3- t^4+ 12t^5- t^6 - 6t^7+t^{10}       \)}}\\
&\multicolumn{6}{l}{\footnotesize{Hilbert Coefficients: $ [ 1, 6, 23, 62, 134, 250, 421, 658, 972, 1374, \ldots ]$}}\\

\evnrow
4 & $ \begin{array}{@{}l@{}} X_{2,3^4,4^4}\\\quad\subset\bP(1^7,3)\end{array}$ & $\footnotesize{\begin{pmatrix}
1 &2 &2\\
1 &2 &2\\
2 &2 &3
\end{pmatrix}}$  & 1/3 & 
      924\newline
    $Y_{2,3^4}
   \subset\bP^6$
 & $\footnotesize{\begin{pmatrix}
1 &1 &1 &2\\
        &1 &1 &2\\
       &&1 &2\\
               &&&2
\end{pmatrix}}$ & $\begin{matrix}
    T_1&10\\
    \textbf{T$_5$} & 8\\
    J_{12} &12\\
    J_{15}&11
\end{matrix}$\\
\evnrow&\multicolumn{6}{l}{\footnotesize{Hilbert Numerator: \( 1- t^2-4t^3+8t^5-4t^7-t^8+t^{10} \)}}\\
 \evnrow&\multicolumn{6}{l}{\footnotesize{Hilbert Coefficients: $ [1, 7, 27, 74, 161, 301, 508, 795, 1175, 1662,\ldots ]$}}\\

 
5 & $ \begin{array}{@{}l@{}} X_{2,3^4,4^4}\\\quad\subset\bP(1^6,2^2)\end{array}$ & $\footnotesize{\begin{pmatrix}
1 &1 &2\\
1 &1 &2\\
2 &2 &3
\end{pmatrix}}$ \\
&\multicolumn{6}{l}{\footnotesize{Hilbert Numerator: \(1- t^2- 4t^3+ 8t^5- 4t^7 - t^8+t^{10}       \)}}\\
&\multicolumn{6}{l}{\footnotesize{Hilbert Coefficients: $ [ 1, 6, 22, 58, 124, 230, 386, 602, 888, 1254, \ldots ]\longleftrightarrow P_{925}(t)$ in \cite{GRDB}.}} \\

\evnrow
6 & $ \begin{array}{@{}l@{}} X_{4^9}\\\quad\subset\bP(1^4,2^4)\end{array}$ & $\footnotesize{\begin{pmatrix}
2 &2 &2\\
2 &2 &2\\
2 &2 &2
\end{pmatrix}}$ &&&&\\
\evnrow&\multicolumn{6}{l}{\footnotesize{Hilbert Numerator: \( 1- 9t^4+16t^6-9t^8+t^{12} \)}}\\
 \evnrow&\multicolumn{6}{l}{\footnotesize{Hilbert Coefficients: $ [1, 4, 14, 36, 74, 140, 234, 364, 536, 756,\ldots ]$}}\\

7 & $\begin{array}{@{}l@{}} X_{3^2,4^5,5^2}\\\quad\subset\bP(1^5,2^2,3)\end{array}$ & \footnotesize{$\begin{pmatrix} 
1 &1 &2\\
2 & 2 &3\\
2 &2 &3
\end{pmatrix}$}  & 1/3 &  928 \newline $Y_{3^2,4^3}\subset \bP(1^5,2^2)$ &
\footnotesize{$\begin{pmatrix}
1 &1 &2 &2\\
  &1 &2 &2\\
     &&2 &2\\
        &&&3
\end{pmatrix}$} & $\begin{matrix}
    T_1 &12\\
\textbf{T$_4$} &10\\
J_{12}&16\\
J_{14}&14\\
J_{45}&13
\end{matrix}$    
\\
&\multicolumn{6}{l}{\footnotesize{Hilbert Numerator: \(1- 2t^3- 5t^4+2t^5+8t^6+2t^7-5t^9-2t^9+t^{12}       \)}}\\
&\multicolumn{6}{l}{\footnotesize{Hilbert Coefficients: $ [ 1, 5, 17, 44, 93, 141, 286, 445, 655, 924, \ldots ]$}}\\

\evnrow
8 & $ \begin{array}{@{}l@{}} X_{4^6,6^3}\\\quad\subset\bP(1^5,3^3)\end{array}$ & $\footnotesize{\begin{pmatrix}
1 &1 &1\\
3 &3 &3\\
3 &3 &3
\end{pmatrix}}$ & 1/3 & 744 \newline $Y_{4^4,6}\subset\bP(1^5,3^2)$ & \footnotesize{$\begin{pmatrix}
1 &1 &1 &1\\
        &3 &3 &3\\
                &&3 &3\\
                        &&&3
\end{pmatrix}$} & $\begin{matrix}
    T_1 &14\\
    \textbf{T$_2$}& 12\\
    J_{23} & 15
\end{matrix}$\\
\evnrow&\multicolumn{6}{l}{\footnotesize{Hilbert Numerator: \( 1- 6t^4+2t^5-3t^6+12t^7-3t^8+2t^9-6t^{10}+t^{14} \)}}\\
 \evnrow&\multicolumn{6}{l}{\footnotesize{Hilbert Coefficients: $ [1, 5, 15, 38, 79, 143, 238, 369, 541, 762,\ldots ]$}}\\


9 & $ \begin{array}{@{}l@{}} X_{2,3^2,5^2,6^4}\\\quad\subset\bP(1^5,2^2,5)\end{array}$ & $\footnotesize{\begin{pmatrix}
1 &1 &2\\
1 &1 &2\\
4 &4 &5
\end{pmatrix}}$  & 1/5 & $Y_{2,3^2,5^2}\newline\subset\bP(1^5,2^2)$ & \footnotesize{$\begin{pmatrix}
-1 &1 &1 &2\\
        &1 &1 &2\\
                &&3 &4\\
                        &&&4
\end{pmatrix}$} &  $\begin{matrix}
    \textbf{T$_5$}&5
\end{matrix}$\\
&\multicolumn{6}{l}{\footnotesize{Hilbert Numerator: \(1- t^2- 2t^3+ 2t^4- 2t^5 - 2t^6+8t^7-2t^8-2t^9+2t^{10}-2t^{11}-t^{12}+t^{14}       \)}}\\
&\multicolumn{6}{l}{\footnotesize{Hilbert Coefficients: $ [ 1, 5, 16, 38, 78, 141, 233, 360, 527, 741, \ldots ]$}}\\
 
\evnrow
10 &  $ \begin{array}{@{}l@{}} X_{3,4^3,5^3,6^2}\\\quad\subset\bP(1^4,2^2,3^2)\end{array}$ & $\footnotesize{\begin{pmatrix}
1 &1 &2\\
2 &2 &3\\
3 &3 &4
\end{pmatrix}}$  & 1/3 & $Y_{3,4^2,5,6}\newline\subset\bP(1^4,2^2,3)$ & $\footnotesize{\begin{pmatrix}
0 &1 &1 &2\\
        & 2&2 &3\\
                &&3 &4\\
                        &&&4
\end{pmatrix}}$ & $\begin{matrix}
    \textbf{T$_2$} &16\\
    J_{23}&18\\
    N&14
\end{matrix}$\\
\evnrow&\multicolumn{6}{l}{\footnotesize{Hilbert Numerator: \( 1- t^3-3t^4-2t^5+2t^6+6t^7+2t^8-2t^9-3t^{10}-t^{11}+t^{14} \)}}\\
 \evnrow&\multicolumn{6}{l}{\footnotesize{Hilbert Coefficients: $ [1, 4, 12, 29, 59, 106, 175, 270, 395, 555,\ldots ]\longleftrightarrow P_{745}(t)$ in \cite{GRDB}}}\\

11 & $\begin{array}{@{}l@{}} X_{2,4^4,6^4}\\\quad\subset\bP(1^5,3^3)\end{array}$ & $\footnotesize{\begin{pmatrix}
1 &1 &3\\
1 &1 &3\\
3 &3 &5
\end{pmatrix}}$ &  1/3 & $Y_{2,4^2,6^2}\newline\subset\bP(1^5,3^2)$ & $\footnotesize{\begin{pmatrix}
-1 &1 &1 &3\\
        &1 &1 &3\\
                &&3 &5\\
                        &&&5
\end{pmatrix}}$ & $\begin{matrix}
 \textbf{T$_3$}& 16\\
 J_{15} &19
\end{matrix}$\\
&\multicolumn{6}{l}{\footnotesize{Hilbert Numerator: \(1- t^2- 4t^4+ 4t^5- 4t^6 + 8t^7-4t^8+4t^9-4t^{10}-t^{12}+t^{14}       \)}}\\
&\multicolumn{6}{l}{\footnotesize{Hilbert Coefficients: $ [ 1, 5, 14, 33, 66, 117, 192, 295, 430, 603, \ldots ]$}}\\
\evnrow
12 & $ \begin{array}{@{}l@{}} X_{2,4^4,6^4}\\\quad\subset\bP(1^6,3,5)\end{array}$ & $\footnotesize{\begin{pmatrix}
1 &1 &3\\
1 &1 &3\\
3 &3 &5
\end{pmatrix}}$ &  1/5 & 926 \newline $Y_{2, 4^4}\subset\bP(1^6,3)$ & $\footnotesize{\begin{pmatrix}
1 &1 &1 &3\\
        &1 &1 &3\\
                &&1 &3\\
                        &&&3
\end{pmatrix}}$ &  $\begin{matrix}
    T_1&8\\
    \textbf{T$_5$} &6\\
    J_{12}& 9
\end{matrix}$\\
\evnrow&\multicolumn{6}{l}{\footnotesize{Hilbert Numerator: \( 1- t^2-4t^4+4t^5-4t^6+8t^7-4t^8+4t^9-4t^{10}-t^{12}+t^{14} \)}}\\
 \evnrow&\multicolumn{6}{l}{\footnotesize{Hilbert Coefficients: $ [1, 6, 20, 51, 107, 197,329,511, 752, 1060,\ldots ]$}}\\


13 & $ \begin{array}{@{}l@{}} X_{4^4,5^4,6}\\\quad\subset\bP(1^4,2^2,3^2)\end{array}$ & $\footnotesize{\begin{pmatrix}
1 &2 &2\\
2 &3 &3\\
2 &3 &3
\end{pmatrix}}$ & 1/3 & 745\newline $Y_{4^3,5^2}\newline\subset\bP(1^4,2^2,3)$ & $\footnotesize{\begin{pmatrix}
1 &2 &2 &2\\
  &2 &2 &2\\
         &&3 &3\\
                &&&3
\end{pmatrix}}$ & $\begin{matrix}
     \textbf{T$_1$} &14\\
     T_3 &12\\
     J_{13}&18\\
     J_{34}&16
\end{matrix}$\\
&\multicolumn{6}{l}{\footnotesize{Hilbert Numerator: \(1- 4t^4- 4t^5+ 3t^6+8t^7+3t^8-4t^9-4t^{10}+t^{14]}       \)}}\\
&\multicolumn{6}{l}{\footnotesize{Hilbert Coefficients: $ [ 1, 4, 12, 30,62, 112, 186, 288, 422, 594, \ldots ]$}}\\


\evnrow
14 &$ \begin{array}{@{}l@{}} X_{3,4^2,5,6^3,7^2}\\\quad\subset\bP(1^4,2^2,3,5)\end{array}$ & $\footnotesize{\begin{pmatrix}
1 &1 &2\\
2 &2 &3\\
4 &4 &5
\end{pmatrix}}$   &1/5 &$Y_{3,4^2,5,6}\newline\subset \bP(1^4,2^2,3)$ &
$\footnotesize{\begin{pmatrix}
0 &1 &1 &2\\
  &2 &2 &3\\
         &&3 &4\\
                &&&4
\end{pmatrix}}$&\\
\evnrow&\multicolumn{6}{l}{\footnotesize{Hilbert Numerator: \( 1- t^3-2t^4-2t^6+t^7+6t^8+t^9-2t^{10}-2t^{12}-t^{13}+t^{16} \)}}\\
 \evnrow&\multicolumn{6}{l}{\footnotesize{Hilbert Coefficients: $ [1, 4, 12, 28, 56, 101, 166, 255, 373, 523,\ldots ]$}}\\


15 & $ \begin{array}{@{}l@{}} X_{4,5^4,6^4}\\\quad\subset\bP(1^3,2^2,3^3)\end{array}$ & $\footnotesize{\begin{pmatrix}
2 &2 &3\\
2 &2 &3\\
3 &3 &4
\end{pmatrix}}$  & 1/3 & 750\newline $Y_{4,5^2,6^2}\newline\subset\bP(1^3,2^2,3^2)$ & $\footnotesize{\begin{pmatrix}
1 &2 &2 &3\\
        &2 &2 &3\\
                &&3 &4\\
                        &&&4
\end{pmatrix}}$ & $\begin{matrix}
    T_1 &18\\
     \textbf{T$_3$} &16\\
     J_{15}& 22\\
     J_{35} &20
\end{matrix}$\\
&\multicolumn{6}{l}{\footnotesize{Hilbert Numerator: \(1- t^4- 4t^5- 4t^6+4t^7 +8t^8+4t^9-4t^{10}-4t^{11}-t^{12}+t^{16}       \)}}\\
&\multicolumn{6}{l}{\footnotesize{Hilbert Coefficients: $ [ 1, 3, 8, 19, 38, 67,110, 169, 246, 345, \ldots ]$}}\\


\evnrow
16 & $ \begin{array}{@{}l@{}} X_{6^9}\\\quad\subset\bP(1^3,3^5)\end{array}$  & $\footnotesize{\begin{pmatrix}
3 &3 &3\\
3 &3 &3\\
3 &3 &3
\end{pmatrix}}$ &   1/3 & 756\newline $Y_{6^5}\subset\bP(1^3,3^4)$ & \footnotesize{$\begin{pmatrix}
3 &3 &3 &3\\
   &3 &3 &3\\
      &&3 &3\\
          &&&3\\
\end{pmatrix}$}& $\begin{matrix}
    \textbf{T$_1$}&18\\
    J_{12} &27
\end{matrix}$\\
\evnrow&\multicolumn{6}{l}{\footnotesize{Hilbert Numerator: \( 1- 9t^6+16t^9-9t^{12}+t^{18} \)}}\\
 \evnrow&\multicolumn{6}{l}{\footnotesize{Hilbert Coefficients: $ [1,3,6, 15, 30, 51, 84, 129, 186, 261,\ldots ]$}}\\

17 & $ \begin{array}{@{}l@{}} X_{2,3^2,7^2,8^4}\\\quad\subset\bP(1^4,2^2,3,7)\end{array}$ & \footnotesize{$\begin{pmatrix}
1 &1 &2\\
1 &1 &2\\
6 &6 &7
\end{pmatrix}$}   &1/7 &$Y_{2,3^2,7^2}\newline\subset\bP(1^4,2^2,3)$  & \footnotesize{$\begin{pmatrix}
    -3 &1 &1 &2\\
    &1 &1 &2\\
    &&5 &6\\
    &&&6
\end{pmatrix}$}\\
&\multicolumn{6}{l}{\footnotesize{Hilbert Numerator: \(1- t^2- 2t^3+ 2t^4- 2t^7 - 2t^8+8t^9-2t^{10}-2t^{11}+2t^{14}-2t^{15}-t^{16}+t^{18}       \)}}\\
&\multicolumn{6}{l}{\footnotesize{Hilbert Coefficients: $ [ 1, 4, 11, 23, 44, 75, 121, 183, 264, 368, \ldots ]$}}\\

\evnrow
18 &$ \begin{array}{@{}l@{}} X_{4^2,6^5,8^2}\\\quad\subset\bP(1^4,3^3,5)\end{array}$ & $\footnotesize{\begin{pmatrix}
1 &1 &3\\
3 &3 &5\\
3 &3 &5
\end{pmatrix}}$ &  1/3 & 753\newline $Y_{4,6^3,8}\newline \subset\bP(1^4,3^2,5)$ & $\footnotesize{\begin{pmatrix}
1 &1 &1 &3\\
        &3 &3 &5\\
                &&3 &5\\
                        &&&5
\end{pmatrix}}$ & $\begin{matrix}
    T_1&20\\
     \textbf{T$_2$} &18\\
     J_{23} &23\\
     J_{25}&21
\end{matrix}$\\
\evnrow
&&& 1/5 & 748 \newline $Y_{4^2,6^3} \subset\bP(1^4,3^3)$ & $\footnotesize{\begin{pmatrix}
1 &1 &3 &3\\
        &1 &3 &3\\
                &&3 &3\\
                        &&&5
\end{pmatrix}}$ & $\begin{matrix}
    T_1&10\\
    \textbf{T$_4$} &8\\
    J_{12}& 13\\
    J_{14}&11 
\end{matrix}$\\
\evnrow&\multicolumn{6}{l}{\footnotesize{Hilbert Numerator: \( 1- 2t^4-5t^6+4t^7-2t^8+8t^9-2t^{10}+4t^{11}-5t^{12}-2t^{14}+t^{18} \)}}\\
 \evnrow&\multicolumn{6}{l}{\footnotesize{Hilbert Coefficients: $ [1, 4, 10, 23, 45, 79, 129, 197, 286, 400,\ldots ]$}}\\


\oddrow
19 & $ \begin{array}{@{}l@{}} X_{2,5^4,8^4}\\\quad\subset\bP(1^4,2^2,3,7)\end{array}$& $\footnotesize{\begin{pmatrix}
1 &1 &4\\
1 &1 &4\\
4 &4 &7
\end{pmatrix}}$ &   bad  sings &&&\\
\oddrow&\multicolumn{6}{l}{\footnotesize{Hilbert Numerator: \( 1- t^2-4t^5+4t^6-4t^8+8t^9-4t^{10}+4t^{12}-4t^{13}-t^{16}+t^{18} \)}}\\
 \oddrow&\multicolumn{6}{l}{\footnotesize{Hilbert Coefficients: $ [1, 4, 11, 25, 50, 87, 143, 219, 318, 446,\ldots ]$}}\\


\evnrow20 & $ \begin{array}{@{}l@{}} X_{2,5^4,8^4}\\\quad\subset\bP(1^5,2,4,7)\end{array}$ & $\footnotesize{\begin{pmatrix}
1 &1 &4\\
1 &1 &4\\
4 &4 &7
\end{pmatrix}}$ & 1/7 & $Y_{2,5^4}\subset\bP(1^5,2,4)$ & $\footnotesize{\begin{pmatrix}
1 &1 &1 &4\\
         &1 &1 &4\\
                &&1 &4\\
                        &&&4
\end{pmatrix}}$ & $\begin{matrix}
    \textbf{T$_5$}&4\\
    J_{12} &6
\end{matrix}$\\
\evnrow&\multicolumn{6}{l}{\footnotesize{Hilbert Numerator: \(1- t^2- 4t^5+ 4t^6-4t^8 +8t^9-4t^{10}+4t^{12}-4t^{13}-t^{16}+t^{18}       \)}}\\
\evnrow&\multicolumn{6}{l}{\footnotesize{Hilbert Coefficients: $ [ 1, 5, 15, 35,71, 127,209, 322, 471, 661, \ldots ]$}}\\


21 & $ \begin{array}{@{}l@{}} X_{4,5^2,6^3,7^2,8}\\\quad\subset\bP(1^3,2^2,3^2,5)\end{array}$ & $\footnotesize{\begin{pmatrix}
1 &2 &3\\
2 &3 &4\\
3 &4 &5
\end{pmatrix}}$ & 1/5 & 750\newline $Y_{4,5^2,6^2}\newline\subset\bP(1^3,2^2,3^2)$ & $\footnotesize{\begin{pmatrix}
1 &2 &2 &3\\
        &2 &2 &3\\
                &&3 &4\\
                        &&&4
\end{pmatrix}}$ & $\begin{matrix}
    T_1 &10\\
    \textbf{T$_5$} &8\\
    J_{13}&12\\
    J_{34}&11
\end{matrix}$\\
&\multicolumn{6}{l}{\footnotesize{Hilbert Numerator: \( 1- t^4-2t^5-3t^6+3t^8+4t^9+3t^{10}-3t^{12}-2t^{13}-t^{14}+t^{18} \)}}\\
 &\multicolumn{6}{l}{\footnotesize{Hilbert Coefficients: $ [1, 3,8, 18, 35,62, 101, 154, 224, 313,\ldots ]$}}\\


\evnrow22  & $ \begin{array}{@{}l@{}} X_{4,5^2,6^3,7^2,8}\\\quad\subset\bP(1^3,2,3^3,4)\end{array}$& $\footnotesize{\begin{pmatrix}
1 &2 &3\\
2 &3 &4\\
3 &4 &5
\end{pmatrix}}$  & 1/3 & $Y_{4,5,6,7,8}\newline\subset\bP(1^3,2,3^2,4)$ & $\footnotesize{\begin{pmatrix}
0 &1 &2 &3\\
        &2 &3 &4\\
                &&4 &5\\
                        &&&6
\end{pmatrix}}$ & $\begin{matrix}
    T_1&24\\
    \textbf{T$_3$} & 22\\
    J_{25}&26\\
    N &30
\end{matrix}$\\
\evnrow&\multicolumn{6}{l}{\footnotesize{Hilbert Numerator: \(1- t^4- 2t^5-3t^6+3t^8 +4t^9+3t^{10}-3t^{12}-2t^{13}-t^{14}+t^{18}       \)}}\\
\evnrow&\multicolumn{6}{l}{\footnotesize{Hilbert Coefficients: $ [ 1, 3, 7, 16, 31, 53,86, 131, 189, 264, \ldots ]\longleftrightarrow P_{752}(t)$ in \cite{GRDB} }}\\


23 & $\begin{array}{@{}l@{}} X_{3,4^2,7,8^3,9^2}\\\quad\subset\bP(1^3,2^2,3^2,7)\end{array}$  & $\footnotesize{\begin{pmatrix}
1 &1 &2\\
2 &2 &3\\
6 &6 &7
\end{pmatrix}}$ & 1/7&$Y_{3,4^2,7,8}\newline\subset\bP(1^3,2^2,3^2)$ &
$\footnotesize{\begin{pmatrix}
-2 &1 &1 &2\\
        &2 &2 &3\\
                &&5 &6\\
                        &&&6
\end{pmatrix}}$&\\
&\multicolumn{6}{l}{\footnotesize{Hilbert Numerator: \( 1-t^3-2t^4+t^5+t^6-t^7-3t^8+t^9+6t^{10}+t^{11}+\cdots+t^{20} \)}}\\
 &\multicolumn{6}{l}{\footnotesize{Hilbert Coefficients: $ [1, 3, 8, 17, 31, 53, 85, 128, 184, 256,\ldots ]$}}\\


\evnrow24 & $ \begin{array}{@{}l@{}} X_{4,5,6^2,7^2,8^2,9}\\\quad\subset\bP(1^3,2,3^2,4,5)\end{array}$ & $\footnotesize{\begin{pmatrix}
1 &2 &3\\
2 &3 &4\\
4 &5 &6
\end{pmatrix}}$ &1/3 & $Y_{4,6,7,8,9}\newline\subset\bP(1^3,2,3,4,5)$ & $\footnotesize{\begin{pmatrix}
0 &1 &2 &4\\
        &2 &3 &5\\
                &&4 &6\\
                        &&&7
\end{pmatrix}}$ & $\begin{matrix}
    \textbf{T$_3$} &26\\
    J_{25} &30\\
    N&28
\end{matrix}$\\

\evnrow
&&& 1/5 & $Y_{4,5,6,7,8}\newline\subset\bP(1^3,2,3^2,4)$ & $\footnotesize{\begin{pmatrix}
0 &1 &2 &3\\
        &2 &3 &4\\
                &&4 &5\\
                        &&&6
\end{pmatrix}}$ &  $\begin{matrix}
    \textbf{T{$_4$}} &10\\
    J_{15}&14\\
    N&12
\end{matrix}$
  \\
\evnrow&\multicolumn{6}{l}{\footnotesize{Hilbert Numerator: \(1- t^4- t^5-2t^6-t^7 +2t^9+4t^{10}+2t^{11}-t^{13}-2t^{14}-t^{15}-t^{16}+t^{20}       \)}}\\
\evnrow&\multicolumn{6}{l}{\footnotesize{Hilbert Coefficients: $ [ 1, 3, 7, 15, 28, 48,77, 116, 167, 232, \ldots ]\longleftrightarrow P_{757}(t)$ in \cite{GRDB}.}}
\label{main table}
\end{longtable}
\end{center}

\bibliographystyle{amsalpha}
\bibliography{References}

\end{document}